\documentclass[10pt]{article}
\usepackage{amstext,amsthm,amsfonts}
\usepackage{amsmath}
\usepackage{geometry}
\usepackage{enumerate}
\usepackage{faktor}
\usepackage[english]{babel}
\usepackage{array}
\usepackage[utf8]{inputenc}
\usepackage{amssymb}
\usepackage{color}
\usepackage{faktor}
\usepackage{xfrac}    
\usepackage{nicefrac}
\usepackage{authblk}

\usepackage{amsopn}
\usepackage{latexsym}
\usepackage{units}
\usepackage{comment}
\usepackage{mathrsfs}

\usepackage{xy}
\usepackage{relsize} 

\def\Si{\Sigma}

\def\lb{\lambda}
\def\Op{\mathfrak{Op}}

\def\Hi{\mathcal{H}}
\def\de{\mathrm{d}}

\def\AA{\mathfrak A}
\def\A{{\mathcal A}}

\def\R{\mathbb R}
\def\C{\mathbb C}

\newtheorem{Theorem}{Theorem}[section]
\newtheorem{proposition}[Theorem]{Proposition}
\newtheorem{Remark}[Theorem]{Remark}
\newtheorem{Definition}[Theorem]{Definition}
\newtheorem{Lemma}[Theorem]{Lemma}
\newtheorem{Corollary}[Theorem]{Corollary}

\author[1]{Fabián Belmonte}
\author[2]{Harold Bustos}
\author[3]{Sebastián Cuellar}
\affil[1,3]{Universidad Católica del Norte}
\affil[2]{Universidad Austral de Chile}
\affil[1]{fbelmonte@ucn.cl}
\affil[2]{harold.bustos@uach.cl}
\affil[3]{sebastian.cuellar01@ucn.cl}
\title{Smooth Fields of Operators and Some Examples Coming from Canonical Quantization}
\date{}

\begin{document}






\maketitle

\begin{abstract}
We introduce a notion of smooth fields of operators following the notion of smooth fields of Hilbert spaces recently defined by L. Lempert and R. Sz\H{o}oke \cite{LS}. Formally, if $\nabla$ is the connection of a smooth field of Hilbert spaces we show that $\hat\nabla=[\nabla,\cdot]$ defines a connection on a suitable space of fields of operators. In order to provide examples we prove that, if $u$ is a suitable constant of motion of $h(q,p)=\|q\|^2$ (i.e.\ $\{u,h\}=0$), then $\Op(u)$ is a smooth field of operators over the open interval $(0,\infty)$, where $\Op$ denotes the canonical quantization (Weyl calculus). Moreover, in such case we show that we can compute derivatives using the formula $\hat\nabla_{X_0}(\Op(u))=\Op(\tilde\nabla_{X_0}(u))$, where $\tilde\nabla$ is a Poisson connection on the Poisson algebra of constants of motion and $X_0=2\lambda\frac{\partial}{\partial \lambda}$. We also introduce a notion of smooth field of $C^*$-algebras and we give an example using Hilbert modules theory.  
  
\end{abstract}

\section{Introduction.}

Let $p:H\to \Lambda$ be a field of Hilbert spaces, i.e. $p$ is a surjective map such that $\Hi(\lb):=p^{-1}(\lb)$ is a Hilbert space, and denote by $\langle\cdot,\cdot\rangle_\lb$ the corresponding inner product. We denote  by $\Gamma$ the set of all sections of such field. For any pair of sections $\varphi,\psi\in\Gamma$, we set the function
$$
h(\varphi,\psi)(\lb)=\langle\varphi(\lb),\psi(\lb)\rangle_\lb.
$$

In order to obtain an interesting mathematical object, we should add further assumptions on a given field of Hilbert spaces. For instance, the notions of measurable and continuous field of Hilbert spaces were introduced by von Neumann \cite{Ne} and Godement \cite{God} respectively (see the appendix \ref{MyC} for details), and they were successfully studied and applied since then. In this article, we are going to introduce a notion of smooth fields of operators following the notion of smooth fields of Hilbert spaces defined in \cite{LS}, and we will provide a large class of examples coming from (Weyl) canonical quantization (including a formula to compute the derivative of suitable fields of operators). Following our results, we will also introduce a notion of smooth field of $C^*$-algebras and we will give an example using Hilbert $C^*$-modules theory.

Let us provide some motivation for our construction. The notions of measurable and continuous field of Hilbert spaces are defined in terms of some space of sections $\Gamma^0$ satisfying certain properties (see definitions \ref{mfhs} and \ref{CField}). However, both notions can be described in terms of bundles or  in terms of Hilbert $C^*$-modules. For instance, if $\Lambda$ is a locally compact Hausdorff space and $p:H\to \Lambda$ is a continuous field of Hilbert spaces, then there is a unique suitable topology on $H$ such that $p:H\to\Lambda$ is a continuous Hilbert bundle and $\Gamma^0$ is contained in the space of continuous sections $\Gamma(\Lambda)$ of $p:H\to\Lambda$. Moreover, the space of continuous sections vanishing at the infinity $\Gamma_0(\Lambda)$ is a Hilbert $C_0(\Lambda)$-module, and conversely, every Hilbert $C_0(\Lambda)$-module can be obtained in this way. There are analogous results in the measurable framework. In appendix \ref{MyC} we explain in more detail those equivalences in both cases.

A similar situation occurs in the $C^*$-algebraic framework. Indeed, upper semi-continuous fields of $C^*$-algebras over $\Lambda$ are in one to one correspondence with upper semi-continuous $C^*$-bundles over $\Lambda$ and with $C_0(\Lambda)$-algebras, for any locally compact Hausdorff space $\Lambda$ (for instance, see theorem C.26 in \cite{Wi}). 

In the smooth framework, we are aware of two concepts so far: the notion of smooth field of Hilbert spaces recently defined in \cite{LS} and the notion of Hermitian bundle (endowed with a connection), but we do not know yet if those notions are equivalent. Let us begin recalling the definition of smooth field of Hilbert spaces and we will explain why we prefer to work with it in this article.

\begin{Definition}\label{def}
Let $\Lambda$ be a finite dimensional smooth manifold and $\operatorname{Vect}(\Lambda)$ the space of smooth vector fields on $\Lambda$. A smooth structure on a field of Hilbert spaces $H\to\Lambda$ is given by specifying a set of sections $\Gamma^\infty$, closed under addition and under multiplication by elements of $C^\infty(\Lambda)$, and a map $\nabla:\operatorname{Vect}(\Lambda)\times\Gamma^\infty\to\Gamma^\infty$ such that, for $X,Y\in\operatorname{Vect}(\Lambda)$, $a\in C^\infty(\Lambda)$ and $\varphi,\psi\in\Gamma^\infty$:
\begin{enumerate}
\item[i)] $\nabla_{X+Y}=\nabla_X +\nabla_Y$, $\nabla_{a X}=a\nabla_X$, $\nabla_X(a\varphi)=X(a)\varphi+a\nabla_X(\varphi)$
\item[ii)] $h(\varphi,\psi)\in C^\infty(\Lambda)$ and $Xh(\varphi,\psi)=h(\nabla_X\varphi,\psi)+h(\varphi,\nabla_{\overline{X}}\psi)$
\item[iii)] $\Hi^\infty(\lb):=\{\varphi(\lb)\mid \varphi\in\Gamma^\infty\}$ is dense in $\Hi(\lb)$, for all $\lb\in\Lambda$.
\end{enumerate} 
\end{Definition}

Conditions \textit{ii)} and \textit{iii)} implies that every smooth field of Hilbert spaces is continuous (see definition \ref{CField}). Moreover, the corresponding space of continuous sections $\Gamma^0(\Lambda)$ is the closure of $\Gamma^\infty$ with the Fréchet topology of local uniform convergence. Similarly, we can complete $\Gamma^\infty$ to obtain the space of $n$-times differentiable sections $\Gamma^n(\Lambda)$ (using the seminorms given by equation \eqref{SN}), for every $n\in\mathbb{N}\cup\{\infty\}$.  

The simplest example of a smooth field of Hilbert spaces is the trivial case, i.e. we take each fiber constant $\Hi(\lb)=V$, $\Gamma^\infty(\Lambda)=C^\infty(\Lambda,V)$ and $\nabla_X=X$. A more interesting case is when $H\to\Lambda$ admits a smooth trivialization.
\begin{Definition}\label{triv1}
Let $p:H\to \Lambda$ be a smooth field of Hilbert spaces with connection $\nabla$. We say that $p:H\to \Lambda$ admits a smooth trivialization if there is a Hilbert space $V$ and a map $T: H\to V$ such that $T_\lb=T |_{\Hi(\lb)}:\Hi(\lb)\to V$ is unitary, $T(\Gamma^\infty)\subseteq C^\infty(\Lambda,V)$ and for each $X\in \operatorname{Vect}(\Lambda)$ we have that
$$
T\nabla_X\varphi=XT\varphi
$$
We say that $T$ is full if in addition we have that $T(\Gamma^\infty(\Lambda))=C^\infty(\Lambda,V)$.
\end{Definition}
Since each $T_\lb$ is unitary, $T$ is continuous with respect to the locally uniform convergence topology, hence $T(\Gamma^n(\Lambda))\subseteq C^n(\Lambda,V)$, for every $n\in\mathbb{N}\cup\{\infty\}$. Therefore, $T$ is full if and only if $T(\Gamma^\infty)$ is dense in $C^\infty(\Lambda,V)$, and in such case we have that $T(\Gamma^n(\Lambda))= C^n(\Lambda,V)$, for every $n\in\mathbb{N}\cup\{\infty\}$.

More generally, we may suppose that the map $T$ locally satisfies the identity 
\begin{equation}\label{pojec}
    T\nabla_X\varphi=XT\varphi+\alpha(X)T\varphi,
\end{equation}
where $\alpha:\operatorname{Vect}(\Lambda)\to \text{End}(V)$ is a suitable map. For instance, if $\alpha$ is a $1$-form (in particular, $\alpha(X)\in C^\infty(\Lambda)$), the couple $(T,V)$ is called a projective trivialization (see definition 2.4.1 in \cite{LS}). It turns out that, every projective trivialization can be transformed into a smooth trivialization by taking the tensor product of the initial smooth field of Hilbert spaces with a suitable line bundle (see subsection 2.4 in \cite{LS}).

The other  concept that we should consider is the notion of Hermitian bundle $p:H\to \Lambda$ \cite{LS} with an Hermitian connection $\nabla$ (allowing the total space to be a Banach manifold and the fibers to be infinite dimensional Hilbert spaces, see subsection 2.1 in \cite{LS} for details). Clearly, the corresponding space of smooth sections $\Gamma^\infty(\Lambda,H)$ defines a smooth field of Hilbert spaces according to definition \ref{def}. However, we do not know yet if for every smooth field of Hilbert spaces $p:H\to\Lambda$ with corresponding section space $\Gamma^\infty$, there is an Hermitian bundle structure on $p:H\to\Lambda$ such that $\Gamma^\infty\subseteq\Gamma^\infty(\Lambda,H)$. Moreover, it seems difficult to construct Hermitian bundles directly. Once a field of Hilbert spaces is given, the main difficulty is to construct the required local trivializations. We expect that the study of smooth fields of Hilbert spaces will become a fundamental tool to overcome that
issue (just as in the continuous or measurable framework). For instance, we know that if $p:H\to \Lambda$ is a smooth field of Hilbert spaces and through every point in $H$ there passes a horizontal section then $p:H\to\Lambda$ is trivializable. Moreover, if $p:H\to\Lambda$ is flat and analytic then the latter condition holds locally (see remark \ref{RM21} or lemma 4.2.1, theorems 2.3.2 and 5.1.2 in \cite{LS}).

Our first goal is to show that the notion of smooth field of Hilbert spaces allows us to introduce a notion of smooth fields of operators and also to define a connection $\hat\nabla$ on a suitable space of fields of operators satisfying properties analogous to \textit{i)} and \textit{ii)} in definition \ref{def}. Formally, if such connection is applied on a field of operators $A=\{A(\lb)\}$, where $A(\lb)$ is an operator on $\Hi(\lb)$ (with suitable domain), then we would expect that the Leibniz's identity holds, i.e.
$$
\nabla_X(A\varphi)=\hat\nabla_X(A)\varphi+A\nabla_X(\varphi),\quad\forall\varphi\in\Gamma^\infty.
$$
Hence, the natural definition for such connection is $\hat\nabla(A)=[\nabla_X,A]=\nabla_X A-A\nabla_X$. The first part of section \ref{S1} is meant to construct suitable spaces of fields of operators $A$ where such expression is well-defined, and $\hat\nabla_X(A)$ is once again a field of operators. Precisely, we will prove that $\hat\nabla_X(A)\in\AA^{n-1}$, for any vector field $X$ and $A\in\AA^n$, where the space of fields of operators $\AA^n$ is introduced in definition \ref{defA}, for any $n\in\mathbb N$ (see theorem \ref{desc}). 

Furthermore, every trivialization $(T,V)$ of $p:H\to\Lambda$ induces a sort of trivialization $\hat{T}$ for $\hat\nabla$ on $\AA^n$. Let us be more precise. Essentially, $\hat{T}$ is defined by $\hat{T}A(\lambda)=\hat{A}(\lambda)=T_\lb A(\lb)T_\lb^*$.\ Notice that each $\hat A(\lb)$ is an operator on $V$, but their corresponding domains might be different depending on $\lambda$. In order to overcome that technical difficulty, we will define a common domain $V^\infty$ (see definition \ref{vinf}). It turns out that $\hat{T}$ maps $\AA^n$ into $C^n(\Lambda,L(V^\infty,V)_{\ast-s})$ and $\hat T\hat \nabla_X=X\hat T$, where $L(V^\infty,V)_{\ast-s}$ denotes the space of linear operators from $V^{\infty}$ to $V$ with the $\ast$-strong topology (see theorem \ref{triviau}).

In section \ref{S1}, we will also consider the trivializable case to give an additional motivation for our construction. Notice that the unitary operators  $U_{\lb_1,\lb_2}=T_{\lb_2}^*T_{\lb_1}:\Hi(\lb_1)\to \Hi(\lb_2)$ play the role of the parallel transport of the connection $\nabla$. 
If $A=\{A(\lb)\}$ is a field of operators, we would like to compare $A(\lb)$ and $U_{\lb,\lb_0}^*A(\lb_0)U_{\lb,\lb_0}$, for any $\lb_0,\lb\in\Lambda$ (both operators have domains in $\Hi(\lb)$). We will show that those operators coincide if and only if we have that $\hat\nabla_X(A)=0$, for any vector field $X$ (see theorem \ref{Teo1}).\ In such case we say that $A$ is an horizontal field of operators. We will also give a weak estimation of the difference between $A(\lb)$ and $U_{\lb,\lb_0}^*A(\lb_0)U_{\lb,\lb_0}$ in the general case (see proposition \ref{est}).

In subsection \ref{IMS}, we will analyze fields of operators belonging to $\AA^n$ as single operators acting on the direct integral $\Hi=\int^\oplus_\Lambda \Hi(\lb)\de\eta(\lb)$, where $\de\eta$ is a fix volume form on $\Lambda$. Since the fields of operators that we are considering belong to some $\AA^n$, the natural domain in the direct integral is $\Gamma_2^\infty=\Gamma^\infty\cap\Hi$, and we will assume that $\Gamma^\infty_2$ is dense in $\Hi$. We will obtain a sort of reduction theorem for our setting analogous to the well known result in the measurable framework \cite{Ne,D,Nu}. More precisely, we will prove that, if $A:\Gamma^\infty_2\to\Gamma^0(\Lambda)\cap\Hi$ is an operator such that $A^*(\Gamma^\infty_2)\subset\Gamma^0(\Lambda)\cap\Hi$, then $A$ can be decomposed through the direct integral if and only if $fA=Af$ on $\Gamma_2^\infty$, for every $f\in C_c^\infty(\Lambda)$ (see theorem \ref{desc1}). We will also introduce a space of fields of operators $\AA^n_2$ such that $\AA^n_2\subseteq\AA^n$ and $\hat\nabla_{X_1}\cdots\hat\nabla_{X_k}(A)$ is an operator on $\Hi$ (with domain $\Gamma^\infty_2$), for every $A\in\AA^n_2$ and $X_1,\cdots, X_k\in\text{Vect}(\Lambda)$, with $0\leq k\leq n$ (see definition \ref{smoothdir}). In the trivializable case we will obtain two further remarkable properties that will become useful later. Let $X\in\text{Vect}(\Lambda)$ and denote by $r_t$ the corresponding flow on $\Lambda$. We will show that the operator $-i(\nabla_X+\frac{1}{2}\text{div}(X))$ is selfadjoint on $\Hi$ (see proposition \ref{HamX}). If $(T,V)$ is a trivialization, define $R_t=T^*r_t^*T$. Notice that, if $A\in\AA^n$, then the map $t\to \langle R_t A R_{-t}\varphi(\lb),\psi(\lb)\rangle_\lb$ is  differentiable for every $\lb\in\Lambda$, and 

$$
\frac{\de}{\de t}\langle R_t A R_{-t}\varphi(\lb),\psi(\lb)\rangle_\lb \Big|_{t=s}=\langle R_s\hat\nabla_X (A)R_{-s} \varphi(\lb),\psi(\lb)\rangle_\lb.
$$

We will prove that if $\Gamma^\infty_2$ is invariant by $R_t$ and $A$ belongs $\AA^1$, then for suitable $\varphi,\psi\in\Gamma^\infty_2$, the map $t\to \langle R_t A R_{-t}\varphi,\psi\rangle$ is  differentiable, and the previous identity holds not only pointwise but also weakly (see theorem \ref{pointw}).    

In section \ref{EX}, we will study an important example, where the smooth field of Hilbert spaces is trivializable by construction, and we will consider fields of operators coming from canonical quantization.  Let $Q^2$ be the selfadjoint operator on $L^2(\R^n)$ given by 
$[Q^2\varphi](q)=||q||^2\varphi(q)$. It is well known that the map  $T:L^2(\R^n)\to\int_{(0,\infty)}^\oplus L^2(\mathbb{S}^{n-1})\de\lb$ given by $T \varphi(\lb,z)=2^{-1/2}\lb^{\frac{n-2}{4}}\varphi(\sqrt{\lb}z)$ is a spectral diagonalization of $Q^2$. Moreover, $T$ can be regarded as a trivialization of the field of Hilbert spaces  $\Hi(\lb)=L^2(\mathbb{S}^{n-1}_{\sqrt\lb},\mu_\lb)$, where $\mu_\lb=2^{-1/2}\sqrt{\lb}\eta_\lb$ and $\eta_\lb$ is the canonical measure on the sphere $\mathbb{S}^{n-1}_{\sqrt\lb}$. Furthermore, we can identify $\varphi\in C^\infty(\R^n)$ with the section given by $\varphi(\lb)=\varphi|_{\mathbb{S}^{n-1}_{\sqrt\lb}}$. Defining $\nabla_X=T^{-1}X T$ on $C_c^\infty(\R^n)$, we obtain a smooth structure on the latter field of Hilbert spaces, and by construction $T$ is a full smooth trivialization. We will give an explicit expression of $\nabla_X$ in proposition \ref{con1} applying a direct computation, and we will also obtain another expression in equation \eqref{con2} using a more geometrical argument.

Recall that canonical reduction theory \cite{D,Nu} implies that if $A$ is a self-adjoint operator $A$ strongly commuting with $Q^2$, then there is a field of self-adjoint operators $\{A(\lb)\}$ such that $TA\varphi(\lb)=A(\lb)T\varphi(\lb)$. The self-adjoint operators strongly commuting with $Q^2$ are sometimes called quantum constants of motion of $Q^2$. Abusing of the notation, we will also call constants of motion of $Q^2$ to the operators satisfying the condition ensuring decomposability according to our reduction theorem (see theorem \ref{desc1}).

In subsection \ref{ehfo}, we will consider operators  of the form $\Op(u)$, where $u\in S'(\R^{2n})$ and $\Op$ is the canonical Weyl quantization \cite{Wey,Fol}.  $\Op$ is meant to map classical observables (i.e. smooth functions on $\R^{2n}$) into quantum observables (i.e. selfadjoint operators on $L^2(\R^n)$) in a physically meaningful way. A classical constant of motion of $h(q,p)=\|q\|^2$ is an smooth function $u\in C^\infty(\R^{2n})$ such that $\{h,u\}=0$, where $\{\cdot,\cdot\}$ is the canonical Poisson bracket on $\R^{2n}$. Notice that $\Op(h)=Q^2$. In certain sense, we will prove that Weyl quantization maps classical constant of motion of $h(q,p)=\|q\|^2$ into quantum constant of motion of $Q^2$ (see theorem \ref{wpcm}). The latter result is interesting on its own right and it does not depend on the smooth structure.

We will show that, if $u$ is a classical constant of motion and $\Op(u)\in\AA^1$ then, in certain sense, the derivative $\hat\nabla_X(\Op(u))$ is also of the form $\Op(\hat u)$ and the function $\hat u$ is obtained in the following geometrical way. Let $\phi\in C^\infty (\R^n)$ given by $\phi(q)=\|q\|^2$. Clearly, the map $\de\phi:T(\R^n\setminus\{0\})\to T(0,\infty)$ is onto and at any point $q\in \R^n\setminus\{0\}$ its kernel is $T_q\mathbb{S}_{\|q\|^2}^{n-1}$. Therefore, for each $X\in\text{Vect}(0,\infty)$ there is a unique vector field $\tilde X$ on $\R^n\setminus\{0\}$ such that $\tilde X(q)\in T_q^\perp\mathbb{S}_{\|q\|^2}^{n-1}$ and $\de\phi(\tilde X)=X$. Also, let $\hat X$ be the Hamiltonian lift of $\tilde X$, i.e. $\hat X$ is the Hamiltonian vector field on $\R^{2n}$ corresponding to the smooth function $h_{\tilde X}$ given by $h_{\tilde X}(q,p)=\langle \tilde X(q),p\rangle$. One of our main results is the following derivation formula
$$
\hat\nabla_{X_0}(\Op(u))=\Op(\hat X_0(u)),
$$
where $X_0$ is the vector field on $(0,\infty)$ given by $X_0(\lb)=\lb\frac{\partial}{\partial\lb}$ (see theorem \ref{Xu}). It is straightforward to show that $\hat X_0(u)$ is also a classical constant of motion. The main reason why the previous formula holds at least for $X_0$ is that the flow of $\hat X_0$ is linear ($\hat r_t^0(q,p)=(e^t q,e^{-t}p)$). In fact, the proof of the previous formula essentially follows from the relation between Weyl quantization and the metaplectic representation, and theorem \ref{pointw} (the metaplectic representation is also key in the proof of theorem \ref{wpcm}). Since $\hat\nabla_{aX_0}=a\hat\nabla_{X_0}$ and $X_0$ is non-degenerate, we can explicitly compute $\hat\nabla_{X}(\Op(u))$ for any vector field $X$. Moreover, the right-hand side of the latter formula contains an interesting geometrical object. Indeed, the map $u\to\hat X(u)$ defines an abelian Poisson connection on the Poisson algebra of classical constants of motion (see the discussion at the end of subsection \ref{ehfo}). In other words, our formula asserts that Weyl quantization exchanges a classical connection with a quantum connection at least for $X=X_0$. 

We will obtain analogous results if we replace $T$ by the spectral diagonalization $T\circ \mathcal F$ of the Laplacian $-\Delta$, where $\mathcal F$ is the Fourier transform (see corollary \ref{Lapl}).  

In subsection \ref{AM}, we will use the results in subsection \ref{ehfo} to provide important and explicit examples of horizontal fields of operators (recall that $A\in\AA^1$ is horizontal if $\hat\nabla_X(A)=0$, for every $X\in\text{Vect}(0,\infty)$). Let $l_{i,j}$ be the classical angular momenta coordinates, i.e.  $l_{i,j}(q,p)=q_ip_j-q_jp_i$ with $1\leq i<j\leq n$. Also let $J:\R^{2n}\to\R^{\frac{n-1}{2}}$ be the map given by $J(q,p)=(l_{1,2}(q,p),\cdots,l_{n-1,n}(q,p))$. We will prove that $a\circ J$ is a classical constant constant motion of $h(q,p)=\|q\|^2$ and of $h_X(q,p)=\langle \tilde X(q),p\rangle$ as well, for all $X\in\text{Vect}(0,\infty)$. It is not difficult to show directly the latter result, but we prefer to show it interpreting $J$ as a moment map and applying a more general result (see proposition  \ref{LP}). In particular, corollary \ref{cor310} implies that $\Op(a\circ J)$ is a horizontal field of operators under suitable conditions (see also corollary \ref{amh}). It is well known that the quantum angular momenta coordinates  $L_{i,j}=Q_i\frac{\partial}{\partial q_j}-Q_j\frac{\partial}{\partial q_i}$ and the total angular momentum operator $L^2=\sum L_{i,j}^2$ are horizontal fields of operators, hence our result is a wide generalization of this fact.\\ \\
In section \ref{LB}, we will study locally uniformly bounded fields of operators and we will introduce a notion of smoothness for fields of $C^*$-algebras. We denote by $\AA_c^n$ the space of all maps $A:\Gamma^\infty\to\Gamma^n(\Lambda)$ which are continuous with respect to the Fréchet topology of $\Gamma^n(\Lambda)$. In order to characterize $\AA^0_c$, we recall that $\Gamma_0^0(\Lambda)$ is a Hilbert $C_0(\Lambda)$-module and we consider the corresponding space of adjointable operators $\AA_0^0$ (see proposition \ref{zero}). In proposition \ref{inter}, we describe how $\hat\nabla$ relates consecutive spaces $\AA^n_c$ and $\AA^{n-1}_c$. That result will lead us to introduce the spaces $\AA^n_{lb}$ (definition \ref{def47}) in such a way that $\hat\nabla$ defines a sort of bounded connection on it (see corollary \ref{cor21}). If $(T,V)$ is a full projective trivialization, then $\AA^n _{lb}$ is isomorphic with $C^n_{lb}(\Lambda,B(V)_{\ast-st})$ (proposition \ref{pro410}).

Inspired by the results of sections \ref{S1} and \ref{LB}, in subsection \ref{SFCA} we define the notion of an smooth  field of $C^*$-algebras (see definition \ref{DefC}) but, first we shall recall some facts of continuous fields of $C^*$- algebras. Reinterpreting some known results found in the literature, we will show at proposition \ref{exako} that the space of compact operators $\mathbb{K}(\Gamma_0)$ on a Hilbert $C_0(\Lambda)$-module $\Gamma_0$, corresponding to a continuous field of Hilbert spaces $p:H\to\Lambda$, defines a continuous structure on the field of $C^*$-algebras $A(\lambda)=\mathbb{K}(\mathcal{H}(\Lambda))$, and $\mathbb{K}(\Gamma_0)$ coincides with the corresponding spaces of continuous sections. That result will allow us to show that there exists an smooth structure on the field of compact operators associated to an smooth field of Hilbert spaces (corollary \ref{cor414}). 

Finally in appendix \ref{MyC} we summarize some well-known facts concerning measurable and continuous fields of Hilbert spaces, emphasizing that both notions have three equivalent ways to be introduced.

\section{Smooth Fields of Operators}\label{S1}

Throughout this article, we will omit in the notation the map $p$, and so we will denote by $H\to\Lambda$ a field of Hilbert spaces.

Let us recall the definition of the space $\Gamma^n(\Lambda)$  given in subsection 3.1 in \cite{LS}, for $n\in\mathbb N\cup\{\infty\}$. The space $\Gamma^0(\Lambda)$ is the $C(\Lambda)$-module of those sections  of $H$ that are locally uniform limits of a sequence in $\Gamma^\infty$. The space $\Gamma^1(\Lambda)$ is the $C^1(\Lambda)$-module of those $\varphi\in \Gamma^0(\Lambda)$ for which there is a sequence $\varphi_j\in\Gamma^\infty$ such that $\varphi_j\to \varphi$ locally uniformly, and for every $X\in \text{Vec}(\Lambda)$, the sequence $\nabla_X\varphi_j$ converges locally uniformly. For such $\varphi$, we can define $\nabla_X\varphi=\lim\nabla_X\varphi_j$ (see lemma 3.1.2 in \cite{LS}). The space $\Gamma^n(\Lambda)$ is defined inductively: $\varphi\in \Gamma^n(\Lambda)$ if $\varphi$ and $\nabla_X\varphi$ belongs to $\Gamma^{n-1}(\Lambda)$, for all $X\in\text{Vect}(\Lambda)$. Finally,
$\Gamma^\infty(\Lambda)=\cap\Gamma^n(\Lambda)$.
The spaces $\Gamma^n(\Lambda)$ and $\Gamma^\infty(\Lambda)$ are Fréchet spaces with the seminorms defined by
\begin{equation}\label{SN}
||\varphi||_{C,X_1,\cdots X_m}=\sup\{||\nabla_{X_1}\cdots\nabla_{X_m}\varphi(\lb)||:\lb\in C\},   
\end{equation}

where $C\subseteq\Lambda$ is compact, $X_1,\cdots,X_m\in\operatorname{Vect}(\Lambda)$ and $m\leq n$ (we can take $X_1,\cdots X_m\in\Xi$, where $\Xi\subset \text{Vect}(\Lambda)$ is finite and generates the tangent plane at each $\lb\in\Lambda$). 

The space $\Gamma^0(\Lambda)$ is by construction the space of continuous sections of $H\to\Lambda$ regarded as a continuous Hilbert bundle (recall that every smooth field of Hilbert spaces is continuous). We denote by $\Gamma_0^0(\Lambda)$ the space of continuous sections vanishing at infinity. As we mentioned in the introduction, the map $h:\Gamma_0^0(\Lambda)\times\Gamma_0^0(\Lambda)\to C_0(\Lambda)$ defines a Hilbert $C_0(\Lambda)$-module structure on $\Gamma_0^0(\Lambda)$. The latter fact will become important in section \ref{LB}.

Let $(T,V)$ be a full trivialization of $H\to\Lambda$ (definition \ref{triv1}) with $\Lambda$ connected. The space $V$ can be identified as the space of horizontal sections $V_0$, i.e. the space of sections $\varphi\in\Gamma^1(\Lambda)$ such that $\nabla_X\varphi=0$, for every $X\in\operatorname{Vect}(\Lambda)$. Indeed, if $f\in V$ is considered as a constant section, then $T^*f\in V_0$. Conversely, if $\varphi$ belongs to $V_0$, then $T\varphi$ is constant, i.e.
\begin{equation}\label{hor1}
T_{\lb_1}\varphi(\lb_1)=T_{\lb_0}\varphi(\lb_0),
\end{equation}
for every $\lb_0,\lb_1\in\Lambda$. 
\begin{Remark}\label{RM21}
{\rm
In particular, since each $T_\lb$ is unitary, through every point in $H$ there passes a horizontal section. Conversely, the proof of theorem 2.3.2 in \cite{LS} shows that if $\Lambda$ is connected and $H\to\Lambda$ is a smooth field of Hilbert spaces  such that through every point in $H$ there passes a horizontal section, then $H\to\Lambda$ admits a full trivialization (indeed, we only need to take $V=V_0$ and $T^*_\lb\varphi=\varphi(\lb)$). In fact, theorem 2.3.2 in \cite{LS} asserts that any flat and analytic field of Hilbert spaces $H$ over a connected and simply connected base $\Lambda$ admits a trivialization. The proof consists of showing that flatness and analiticity imply that through every point in $H$ there passes a locally defined horizontal section, and since $\Lambda$ is simply connected, that property holds globally (see lemma 4.1.3 and lemma 4.2.1. in \cite{LS}).

}
\end{Remark}

Let us recall one of our motivation to consider the trivializable case. The map $U_{\lb_1,\lb_2}=T_{\lb_2}^*T_{\lb_1}:\Hi(\lb_1)\to \Hi(\lb_2)$ is unitary and plays the role of the parallel transport of the connection $\nabla$. Notice that $U_{\lb_2,\lb_3}U_{\lb_1,\lb_2}=U_{\lb_1,\lb_3}$. If $A=\{A(\lb)\}$ is a field of operators, we would like to compare $A(\lb)$ and $U_{\lb,\lb_0}^*A(\lb_0)U_{\lb,\lb_0}$, for any $\lb_0,\lb\in\Lambda$ (both operators have domains in $\Hi(\lb)$). The following simple result characterize when those operators coincide.

\begin{Theorem}\label{Teo1}
Let $H\to\Lambda$ be a smooth field of Hilbert spaces with connection $\nabla:{\rm Vect}(\Lambda)\times \Gamma^\infty\to\Gamma^\infty$ and let $T:H\to V$ be a smooth trivialization. Also, let $A=\{A(\lb)\}$ be a field of operators such that $A(\Gamma^\infty(\Lambda))\subseteq\Gamma^1(\Lambda)$. The following statements are equivalent.
\begin{enumerate}
\item[a)] $[\nabla_X,A]=0$, for every $X\in\operatorname{Vect}(\Lambda)$.
\item[b)] $AV_0\subseteq V_0$.
\item[c)] The field of operators $\hat A=\{\hat A(\lb)\}$ on $V$ defined by $\hat A(\lb)f=T_\lb A(\lb)T_\lb^*f$ is constant.
\item[d)] $A(\lb)=U_{\lb,\lb_0}^*A(\lb_0)U_{\lb,\lb_0}$, for any $\lb_0,\lb\in\Lambda$.
\end{enumerate}
\end{Theorem}
\begin{proof}
$a)\Rightarrow b)$: Let $\varphi\in V_0$. Then, $\nabla_X A\varphi=A\nabla_X\varphi=0$, i.e. $A\varphi$ is a horizontal section.\\
$b)\Rightarrow c)$: For each $f\in V$, we have that $T^*f\in V_0$, thus $AT^*f\in V_0$. Therefore, equation \eqref{hor1} implies 
$$
 T_\lb A(\lb) T^*_\lb f=T_{\lb_0} A(\lb_0) T^*_{\lb_0} f
$$
for every $\lb_0,\lb\in\Lambda$.\\
$c)\Rightarrow d)$:
$$
U_{\lb,\lb_0}^*A(\lb_0)U_{\lb,\lb_0}=T_{\lb}^*T_{\lb_0}A(\lb_0)T_{\lb_0}^*T_{\lb}=T_{\lb}^*\hat A(\lb_0)T_{\lb}=T_{\lb}^*\hat A(\lb)T_{\lb}=A(\lb).
$$
$d)\Rightarrow c)$: Since $A(\lb)=T_{\lb}^*T_{\lb_0}A(\lb_0)T_{\lb_0}^*T_{\lb}$, we have that
$$
\hat A(\lb)=T_{\lb}A(\lb)T_{\lb}^*=T_{\lb_0}A(\lb_0)T_{\lb_0}^*=\hat A(\lb_0).
$$
$c)\Rightarrow a)$: By definition, $T(A\varphi)=\hat A T\varphi$, for every $\varphi\in\Gamma^\infty(\Lambda)$. Since we are assuming that $\hat A$ is a constant field of operators, we have that
$$
T(A\nabla_X\varphi)=\hat A T\nabla_X\varphi=\hat A X T\varphi=X\hat A T\varphi=XT(A\varphi)=T(\nabla_X A\varphi),
$$
for every $\varphi\in\Gamma^\infty(\Lambda)$. 
\end{proof}

The condition $A(\Gamma^\infty(\Lambda))\subseteq\Gamma^1(\Lambda)$ in the previous theorem guarantee that $A\varphi$ is well-defined and belongs to $\Gamma^1(\Lambda)$, for every $\varphi\in V_0$, and it also implies that the domain of each $A(\lb)$ is the whole Hilbert space $\Hi(\lb)$. If we only assume that $A(\Gamma^\infty)\subseteq\Gamma^1(\Lambda)$, then we would find a problem with the domains involved. We shall address and solve that problem later (see definition \ref{vinf} and remark \ref{S22}). Another possibility is to impose a topological condition (without assuming that $H\to\Lambda$ admits a smooth trivialization).
\begin{proposition}
Let $A=\{A(\lb)\}$ be a field of operators such that $A(\Gamma^\infty)\subseteq\Gamma^1(\Lambda)$ and $[\nabla_X,A]=0$ in $\Gamma^\infty$, for every $X\in\operatorname{Vect}(\Lambda)$. Assume that $A:\Gamma^\infty\to \Gamma^1(\Lambda)$ is continuous with respect to the locally uniform limit topology (i.e. the topology of $\Gamma^0(\Lambda)$). Then, for each $n\in\mathbb N\cup\{\infty\}$,
\begin{enumerate}
    \item[a)] $A$ extends to each $\Gamma^n(\Lambda)$ and $A(\Gamma^n(\Lambda))\subseteq \Gamma^n(\Lambda)$,
    \item[b)] $[\nabla_X,A]=0$ in $\Gamma^n(\Lambda)$ and
    \item[c)] $A:\Gamma^n(\Lambda)\to \Gamma^n(\Lambda)$ is continuous. 
\end{enumerate}
\end{proposition}
\begin{proof}
The continuity of $A$ with respect to the locally uniform topology implies that $A$ extends to a continuous operator on $\Gamma^0(\Lambda)$. In particular, $A\varphi$ is well-defined for every $\varphi\in \Gamma^n(\Lambda)$. Let us show the case $n=1$. The general case $n>1$ follows by induction and the same argument used for $n=1$.

Let $\varphi\in\Gamma^1(\Lambda)$. Hence, there is $\varphi_j\in\Gamma^\infty$ such that $\varphi_j\to\varphi$ and $\nabla_X\varphi_j\to\nabla_X\varphi$. The continuity of $A$ implies that $A\varphi_j\to A\varphi$ in $\Gamma^0(\Lambda)$. Moreover, $\nabla_X A\varphi_j=A\nabla_X\varphi_j\to A\nabla_X\varphi$ in $\Gamma^0(\Lambda)$. Thus, $A\varphi\in \Gamma^1(\Lambda)$ and $\nabla_X A\varphi=A\nabla_X \varphi$. The latter identity implies c) and the continuity of $A:\Gamma^1(\Lambda)\to\Gamma^1(\Lambda)$.  
\end{proof}

We will show that the commutator $\hat\nabla_X(A)=[\nabla_X,A]$ defines a connection on a suitable space of fields of operators. The main purpose of this article is to study $\hat\nabla$. Theorem \ref{Teo1} heuristically asserts that ``the derivative of a field of operators vanishes iff the field of operators is constant''. Let us introduce the basic sets of fields of operators that we will consider in this article.

\begin{Definition}\label{defA}
Let $H\to\Lambda$ a smooth field of Hilbert spaces with connection $\nabla:\operatorname{Vect}(\Lambda)\times \Gamma^\infty\to\Gamma^\infty$. 
We denote by $\AA^n$ the space formed by the fields of operators $A=\{A(\lb)\}$ such that
\begin{enumerate}
\item[i)] the domain of $A(\lb)$ and $A^*(\lb)$ contains $\Hi^\infty(\lb)=\{\varphi(\lb)\mid \varphi\in\Gamma^\infty\}.$
\item[ii)]  $A(\Gamma^\infty)\subseteq \Gamma^n$ and $A^*(\Gamma^\infty)\subseteq \Gamma^n$.
\end{enumerate}

We say that $A=\{A(\lb)\}$ is smooth
 if $A\in \mathfrak{A}^n$, for every $n\in\mathbb{N}$. We denote by $\mathfrak{A}^\infty$ the space formed by the smooth fields of operators. We say that $A\in \mathfrak{A}^1$ is a horizontal field of operators if $\hat\nabla_X(A)=0$, for every $X\in\operatorname{Vect}(\Lambda)$.  
\end{Definition}

In particular, if $A\in\AA^n$ then each $A(\lb)$ is closable. Moreover, it is clear that $\hat\nabla_X(A)(\Gamma^\infty)\subseteq \Gamma^{n-1}$, but it is not so obvious that $\hat\nabla_X(A)$ is given by a field of operators.
The following result will be useful to show that claim.
\begin{Lemma}\label{sic}
Let $H\to\Lambda$ a smooth field of Hilbert spaces with connection $\nabla$. The following identification of quotient spaces holds
\begin{enumerate}
    \item[a)] $\Gamma^{\infty}/K^{\infty}(\lambda)\cong \Hi^{\infty}(\lambda)$, where $K^{\infty}(\lambda):=\{\varphi\in\Gamma^\infty\mid \varphi(\lambda)=0\}$.
    \item[b)] $\Gamma^{n}(\Lambda)/K^{n}(\lambda)\cong \Hi^{n}(\lambda)$, where $K^{n}(\lambda)=\{\varphi\in\Gamma^n(\Lambda)\mid \varphi(\lambda)=0\}$ and $\Hi^{n}(\lambda):=\{\varphi(\lambda)\mid\varphi\in\Gamma^n(\Lambda)\}$. 
    \item[c)] $\Gamma^{0}_0(\Lambda)/K_{0}(\lambda)\cong \Hi(\lambda)$, where $K_{0}(\lambda):=\{\varphi\in\Gamma^0_0(\Lambda)\mid \varphi(\lambda)=0\}$ and  $\Gamma^{0}_0(\Lambda)/K_{0}(\lambda)$ is endowed with the canonical quotient norm.
\end{enumerate}
Let $A:\Gamma^\infty\to\Gamma^n(\Lambda)$ be a linear operator. There is a field of operators $\{A(\lambda)\}$ such that the domain of $A(\lambda)$ contains $\Hi^{\infty}(\lambda)$ and $A\varphi(\lambda)=A(\lambda)\varphi(\lambda)$, for every $\lambda\in\Lambda$ and $\varphi\in\Gamma^\infty$, if and only if $A(K^\infty(\lambda))\subseteq K^n(\lambda)$, for every $\lambda\in\Lambda$.
\end{Lemma}

\begin{proof}
For each $\lb\in\Lambda$, define $\rho_\lb(\varphi)=\varphi(\lb)$, for every section $\varphi$ of $H\to\Lambda$. The restriction of $\rho_\lb$ defines an epimorphism from $\Gamma^\infty$ onto $\Hi^\infty(\lb)$ and its kernel is $K^\infty(\lb)$, and this implies a). The same argument shows b). Let $\varphi\in\Gamma_0^0(\Lambda)$ and $\lb\in\Lambda$. Clearly,
$$
\|\varphi(\lb)\|\leq\inf\Big\{\sup_{\mu\in\Lambda}\|(\varphi-\psi)(\mu)\| \mid \psi\in K_0(\lb)\Big\}.
$$
In order to show that we actually have an equality, it is enough to prove that, for each $\epsilon>0$, there is $\psi\in K_0(\lb)$ such that $\sup_{\mu\in\Lambda}\|(\varphi-\psi)(\mu)\|<\|\varphi(\lb)\|+\epsilon$. Let $C$ be the compact set given by 
$$
C=\{\mu\in\Lambda\mid\|\varphi(\mu)\|\geq \|\varphi(\lb)\|+\epsilon\}.
$$
If $f\in C_0(\Lambda)$ is such that $0\leq f\leq 1$, $f(\lb)=0$ and $f|_C=1$, then $\psi=f\varphi$ satisfies the required inequality. In particular, the subspace $\Hi_0(\lb)=\{\varphi(\lb)\mid \varphi\in\Gamma_0^0(\Lambda)\}$ is closed in $\Hi(\lb)$. Since $\Hi^\infty(\lb)\subset \Hi_0(\lb)$, condition iii) in definition \ref{def} implies c). The last claim of our lemma is a direct consequence of a) and b).
\end{proof}

\begin{Theorem}\label{desc}
Let $H\to\Lambda$ a smooth field of Hilbert spaces with connection $\nabla$.\ For each $n\geq 1$, the map $\hat\nabla:\operatorname{Vect}(\Lambda)\times \mathfrak{A}^n\to\mathfrak{A}^{n-1}$ given by
$$
\hat\nabla_X(A)=[\nabla_X,A]
$$
is well-defined and satisfies the following properties for all $X,Y\in \operatorname{Vect}(\Lambda)$, $a\in C^{\infty}(\Lambda)$.
\begin{enumerate}
\item[i)] $\hat\nabla_{X+Y}(A)=\hat\nabla_{X}(A)+\hat\nabla_{Y}(A),\qquad \hat\nabla_{aX}(A)=a\hat\nabla_{X}(A)$ 
\item[ii)] $\hat\nabla_{X}(aA)=X(a)A+a\hat\nabla_{X}(A),$
\item[iii)] $h(\hat\nabla_{X}(A)\varphi,\psi)=h(\varphi,\hat\nabla_{\overline{X}}(A^{*})\psi)$, for every $\varphi,\psi\in\Gamma^\infty$.
\item[iv)]$\hat\nabla_{\overline{X}}(A^{*})(\lb)\subseteq[\hat\nabla_{X}(A)(\lb)]^*$, for each $\lb\in\Lambda$.
\end{enumerate}
\end{Theorem}
\begin{proof}
We will show first properties \textit{i),ii)} and \textit{iii)} and later that $\hat\nabla_X(A)$ is given by a field of operators satisfying condition a) in definition \ref{defA}. 
The first equality in \textit{i)} is clear. Since $A$ is given by a field of operators, we have that $aA=Aa$ and 
$$
\hat\nabla_{aX}(A)=[\nabla_{aX},A]=[a\nabla_{X},A]=a\hat\nabla_{X}(A).
$$
A direct computation implies \textit{ii)}. 
For \textit{iii)}, a repeated application of \textit{ii)} in definition 1.1.\ gives
$$
Xh(A\varphi,\psi)=h(\nabla_X A\varphi,\psi)+h(A\varphi,\nabla_{\overline{X}}\psi)
$$
$$
=Xh(\varphi,A^{*}\psi)=h(\nabla_X \varphi,A^{*}\psi)+h(\varphi,\nabla_{\overline{X}}A^{*}\psi),
$$
for any two sections $\varphi,\psi \in \Gamma^{\infty}$. Therefore,
$$
h(A\nabla_X \varphi,\psi)+h(\varphi,\nabla_{\overline{X}}A^{*}\psi)=h(\nabla_X A\varphi,\psi)+h(\varphi,A^{*}\nabla_{\overline{X}}\psi)
$$
$$
h(\hat\nabla_{X}(A)\varphi,\psi)=h(\varphi,\hat\nabla_{\overline{X}}(A^{*})\psi).
$$

The last identity implies that $\hat\nabla_{X}(A)(K^\infty(\lambda))\subseteq K^n(\lambda)$, and therefore lemma \ref{sic} shows that  $\hat{\nabla}_X$ is well-defined. 

Clearly, identity \textit{iii)} also implies \textit{iv)}. 
\end{proof}

In order to obtain further properties of fields of operators belonging to $\AA^n$, we should impose additional conditions. For instance, in section \ref{LB} we will consider locally uniformly bounded fields of operators belonging to $\AA^n$. In subsection \ref{IMS} we analyze  field of operators belonging to $\AA^n$ as a single operator acting on a direct integral. 

Let us assume that our smooth field of Hilbert spaces admits a projective trivialization $T:H\to V$. In such case $\hat A(\lb)=T_\lb A(\lb)T_\lb^*$ is a field of operators on $V$. We would like to fix a common domain for those operators. A first candidate might be the space formed by all the vectors $f\in V$ such that $T^*f\in \Gamma^\infty$. However, that space may be trivial (for instance see the example in section \ref{EX}). Instead, we will consider the space where such condition holds locally. 
\begin{Definition}\label{vinf}
Let $T:H\to V$ be a projective trivialization  of the smooth field of Hilbert spaces $H\to\Lambda$. Also, let $f\in V$. We say that $f\in V^\infty$ if  for every  $\lb_0\in\Lambda$, there are $U\subseteq \Lambda$ open and $\varphi\in\Gamma^\infty$ such that $\lb_0\in U$ and $T^*_\lb f=\varphi(\lb)$, for every $\lb\in U$.
\end{Definition}
For instance, in the example of section \ref{EX}, we have that $V=L^2(\mathbb S^{n-1})$ and $V^\infty=C^\infty(\mathbb S^{n-1})$.

By definition, if $A\in\AA^n$ and $f\in V^\infty$, then $A(T^*f)$ is well-defined and belongs $\Gamma^n(\Lambda)$; therefore $\hat A f\in C^n(\Lambda,V)$. The same happens with $A^*$. In other words, if $A\in\AA^n$ then $\hat A\in C^n(\Lambda,L(V^\infty,V)_{\ast-s})$, where $L(V^\infty,V)_{\ast-s}$ is the space of linear operators from $V^{\infty}$ to $V$ with the $\ast$-strong topology. 
\begin{Theorem}\label{triviau}
Let $T:H\to V$ be a projective trivialization of the smooth field of Hilbert spaces $H\to\Lambda$ with connection $\nabla$. The map $\hat T:\AA^n\to C^n(\Lambda,L(V^\infty,V)_{\ast-s})$, given by
$$
\hat T A(\lb)=\hat A(\lb)=T_\lb A(\lb)T^*_\lb,
$$
defines a trivialization of $\AA^n$ with respect to the connection $\hat\nabla$, i.e. the following identity holds
$$\hat{T}(\hat{\nabla}_XA)=X\hat{T}A.
$$

\end{Theorem}

\begin{proof}
Let $A\in\AA^n$, $X\in\operatorname{Vect}(\Lambda)$ and $f\in V^\infty$. Since $\hat A f\in C^n(\Lambda,V)$ and $Xf=0$, we have that 
\begin{eqnarray*}
\hat{T}(\hat{\nabla}_X A) (\lambda) f &=&
(T\nabla_X AT^*)(\lambda)f-(TA\nabla_X T^*)(\lambda)f
\\
&=& (XTAT^*+\alpha(X)TAT^*)(\lambda)f - (TAT^*)(T\nabla_X T^*)(\lambda)f
\\
&=& X\hat{T}A(\lambda)f+\alpha(X)\hat{T}A(\lambda)f - \hat{T}A(\lambda)(Xf+\alpha(X)(\lambda)f)
\\
&=& X\hat{T}A(\lambda)f+\alpha(X)\hat{T}A(\lambda)f - \hat{T}A\alpha(X)(\lambda)f
\\
&=&X\hat{T}A(\lambda)f \,.
\end{eqnarray*}
\end{proof}
\begin{Remark}\label{S22}
{\rm
Theorem \ref{Teo1} still holds when $A\in\AA^1$, $V_0$ is replaced by $T^*V^\infty$ and the domains of the operators $A(\lb)$ and $\hat A(\lb)$ are $T_\lb^*V^\infty$ and $V^\infty$ respectively. 
}
\end{Remark}
\begin{Remark}
{\rm
Notice that if $T$ is trivialization and $V^\infty$ is dense in $V$ then $T$ is full. Indeed, in such case every constant section $f\in V$ belongs to the closure of $T(\Gamma^\infty)$, and the latter implies that $T(\Gamma^\infty)$ is dense in $C^\infty(\Lambda,V)$.
}
\end{Remark}

Let us return to our initial motivation: comparing $A(\lb)$ and $U_{\lb,\lb_0}^*A(\lb_0)U_{\lb,\lb_0}$, where $U_{\lb,\lb_0}=T_{\lb_0}^*T_\lb$. Assume that $A\in\AA^1$. Thus, the map $\lb\to \langle \hat A(\lb)f,g\rangle_V$ belongs to $C^1(\Lambda)$, for every $f,g\in V^\infty$. Let $X\in\operatorname{Vect}(\Lambda)$ and $r_t$ its one parameter flow (integral curve). Therefore, if $\lb=r_t(\lb_0)$
$$
\langle [A(\lb)-U_{\lb,\lb_0}^*A(\lb_0)U_{\lb,\lb_0}]T^*_\lb f,T^*_\lb g\rangle_\lb=\langle [\hat A(\lb)-\hat A(\lb_0)] f,g\rangle_V=\int_0^t \langle X\hat A(r_s\lb_0)f,g\rangle_V\,\de s. 
$$
The following result is a direct consequence of the previous identity.
\begin{proposition}\label{est}
Let $T:H\to V$ be a trivialization of the smooth field of Hilbert spaces $H\to\Lambda$ with connection $\nabla$. Also, let $X\in\operatorname{Vect}(\Lambda)$ and $r_t$ its one parameter flow. Fix $\lb_0\in\Lambda$ and $\lb_1=r_t(\lb_0)$, for some $t>0$. If $A\in\AA^1$ and $v,w\in T^*_{\lb_1}(V^\infty)$, then
$$
|\langle [A(\lb_1)-U_{\lb_1,\lb_0}^*A(\lb_0)U_{\lb_1,\lb_0}]v,w\rangle_{\lb_1}|\leq t \sup_{\lb\in\gamma[0,t]}|\langle\hat\nabla_XA(\lb)(U_{\lb_{1},\lb}v),U_{\lb_{1},\lb}w\rangle_\lb|,
$$
where $\gamma[0,t]=\{r_s(\lb_0):s\in[0,t]\}$ is the integral curve between $\lb_0$ and $\lb_1$. 
\end{proposition}

Our construction may help us to find general trivializations satisfying equation \eqref{pojec}. Indeed, let us consider the following more general framework. Assume that $H_1\to\Lambda$ and $H_2\to \Lambda$ are two smooth fields of operators with connections $\nabla^1$ and $\nabla^2$ respectively. We say that a field of operators $A=\{A(\lb)\}$ belongs to $\AA^n(H_1,H_2)$ if: 
\begin{enumerate}
\item[a)] the domain of $A(\lb)$ contains $\Hi_1^\infty(\lb)$ and the domain of $A^*(\lb)$ contains $\Hi_2^\infty(\lb)$.
\item[b)]  $A(\Gamma_1^\infty)\subseteq \Gamma_2^n$ and $A^*(\Gamma_2^\infty)\subseteq \Gamma_1^n$.
\end{enumerate}
Clearly, we can adapt the proof of theorem \ref{desc} and define $\hat\nabla:\operatorname{Vect}(\Lambda)\times \AA^n(H_1,H_2)\to\AA^{n-1}(H_1,H_2)$ by
$$
\hat\nabla_X(A)=\nabla^2_X A-A\nabla^1_X.
$$
Let $U=\{U(\lb)\}$ be a smooth field of unitary operators (i.e. $U(\Gamma^\infty_1(\Lambda))\subseteq \Gamma^\infty_2(\Lambda)$) and $T_2:H_2\to V$ be a smooth trivialization of $H_2\to \Lambda$. Therefore, the map $T_1=T_2\circ U:H_1\to V$ satisfies the identity
$$
T_1\nabla^1_X\varphi=X T_1\varphi+\alpha(X) T_1\varphi,
$$
where $\alpha(X)=-T_1U^*\hat\nabla_X(U)T_1^*$. 

Notice that the previous expression generalize the projective case. As we mentioned before, if $H\to \Lambda$ admits a projective trivialization, then we can tensorize $H$ with a suitable line bundle $L$ to obtain  smooth field of Hilbert spaces $H\otimes L\to \Lambda$ admitting a smooth trivialization. The fibers $\Hi'(\lb)$ of $H\otimes L\to \Lambda$ as sets coincide with the fibers $\Hi(\lb)$ of $H\to V$, but the corresponding inner product is obtained by multiplying $\langle\cdot,\cdot\rangle_\lb$ by $a(\lb)$, where $a$ is certain non-negative smooth function on $\Lambda$. In particular, $U_\lb=\sqrt{a(\lb)}I_\lb$ is a smooth field of unitary operators, and we recover the smooth trivialization of $H\otimes L\to\Lambda$ taking $T'=T\circ U$.

\subsection{Direct integrals and decomposable operators}\label{IMS}

Historically, fields of Hilbert spaces and operators emerged during the development of the reduction theory of von Neumann algebras \cite{Ne,D}. One key element on that theory was the notion of direct integral of measurable fields of Hilbert spaces. The canonical reduction theorem characterize bounded operators on the direct integral defined by bounded measurable fields of operators. Such result was extended for (unbounded) closed operators in \cite{Nu}. A similar reduction theory was developed in \cite{Ta} in the continuous framework (but certain unnecessary topological condition was assumed). An interesting approach for the unbounded case was recently described in \cite{ABr}. We will show an analogous reduction theorem within our smooth framework, and we will use it to construct examples of smooth fields of operators in section \ref{EX}.

In order to construct direct integrals, we need to endow $H\to\Lambda$ with a measurable structure, or equivalently, to consider a measurable Hilbert bundle structure on $H\to\Lambda$ (see definition 2.4.8 in \cite{We} and appendix \ref{MyC}). For the purposes of this article, we shall assume the following stronger assumption.

\begin{Definition}
Let $H\to \Lambda$ be a smooth field of Hilbert spaces. A weakly smooth trivialization is a couple $(V,T)$, where $V$ is a Hilbert space and $T:H\to V$ is a map such that $T|_{\Hi(\lb)}$ is unitary for any $\lambda\in\Lambda$ and $T(\Gamma^\infty)\subseteq C^\infty(\Lambda,V)$.
\end{Definition}
Let $\eta$ be a Borel measure on $\Lambda$ (for instance, we can fix a density on $\Lambda$ or a volume form $\de\eta$). Recall that the direct integral $
\Hi= \int_\Lambda^\oplus \Hi(\lb)\de\eta(\lb)$ is the Hilbert space of all measurable sections $\varphi$ such that $\int_{\Lambda}||\varphi(\lb)||_\lb^2\de\eta(\lb)<\infty$ (See appendix \ref{MyC} for details). 

Clearly, $L^\infty(\Lambda, \eta)$ can be represented faithfully in $\Hi$ (acting as constant operators on each fiber) and every measurable essentially uniformly bounded field of operators $A=\{A(\lb)\}$ defines a bounded operator in $\Hi$ given by $A\varphi(\lb):=A(\lb)\varphi(\lb)$. Moreover, such operator belongs to the commutant of $L^\infty(\Lambda,\eta)$.
Conversely, every bounded operator on $\Hi$ belonging to the commutant of $L^\infty(\Lambda,\eta)$ can be decomposed as a measurable essentially uniformly bounded field of operators (for instance, see \cite{D} or \cite{BS}). The  extensions of the later result given in \cite{Nu, Ta} do not guarantee that  the domain  
of each $A(\lb)$ contains $\Hi^\infty(\lb)$. We will overcome that issue in our framework. Let us consider the spaces
$$
\Gamma^\infty_2=\Gamma^\infty\cap\Hi=\Big\{\varphi\in\Gamma^\infty\mid \int_{\Lambda}||\varphi(\lb)||_\lb^2\de\eta(\lb)<\infty\Big\}
$$
$$
\Gamma^0_2(\Lambda)=\Gamma^0(\Lambda)\cap\Hi=\Big\{\varphi\in\Gamma^0(\Lambda)\mid \int_{\Lambda}||\varphi(\lb)||_\lb^2\de\eta(\lb)<\infty\Big\}
$$
Through the rest of this article we will assume that $\Gamma^\infty_2$ is dense in $\Hi$ (for instance, this holds true if $T$ is full). Using a suitable bump function, it is straightforward to show that $\{\varphi(\lb)\mid\varphi\in\Gamma^\infty_2\}=\Hi^\infty(\lb)$ and $\{\varphi(\lb)\mid\varphi\in\Gamma^0_2(\Lambda)\}=\Hi(\lb)$. Therefore $\Gamma^\infty_2/K^\infty_2(\lb)=\Hi^\infty(\lb)$ and $\Gamma^0_2(\Lambda)/K^0_2(\lb)=\Hi(\lb)$, where $K^\infty_2(\lb)=\{\varphi\in\Gamma^\infty_2\mid\varphi(\lb)=0\}$ and $K^0_2(\lb)=\{\varphi\in\Gamma^0_2(\Lambda)\mid\varphi(\lb)=0\}$. 

\begin{Theorem}\label{desc1}
Let $H\to\Lambda$ be a smooth field of Hilbert spaces and $A:\Gamma^\infty_2\to\Gamma^0_2(\Lambda)$ be a linear operator. Assume that the domain of $A^*$ contains $\Gamma^\infty_2$ and $A^*(\Gamma^\infty_2)\subseteq\Gamma^0_2(\Lambda)$. There is a field of operators $\{A(\lb)\}$ such that the domain of $A(\lb)$ contains $\Hi^\infty(\lb)$ and $A\varphi(\lb)=A(\lb)\varphi(\lb)$, for every $\varphi\in \Gamma^\infty_2$ and $\lb\in\Lambda$, if and only if $fA=Af$ on $\Gamma^\infty_2$, for every $f\in C_c^\infty(\Lambda)$. In such case, $A$ belongs to $\AA^0$ and moreover, $A$ is locally uniformly bounded if and only if $fA$ is bounded in $\Hi$, for each $f\in C_c^\infty(\Lambda)$.
\end{Theorem}

\begin{proof}
Obviously if $A\varphi(\lb)=A(\lb)\varphi(\lb)$, then clearly the required commutation holds. Let us prove the converse.
Notice that
$$
\int_\Lambda f h(A\varphi,\psi)\de\eta =\langle fA\varphi,\psi\rangle=\langle f\varphi,A^*\psi\rangle=\int_\Lambda f h(\varphi,A^*\psi)\de\eta. 
$$
Therefore $h(A\varphi,\psi)=h(\varphi,A^*\psi)$ almost everywhere. Since both functions are continuous, the latter identity holds everywhere. Hence, $AK^\infty_2(\lb)\subseteq K_2(\lb)$ and we can repeat the argument of lemma \ref{sic} to show that there is a field of operators $\{A(\lb)\}$ such that the domain of $A(\lb)$ contains $\Hi^\infty(\lb)$ and $A\varphi(\lb)=A(\lb)\varphi(\lb)$, for every $\varphi\in \Gamma^\infty_2$ and $\lb\in\Lambda$. In particular, $A$ is well-defined on $\Gamma^\infty$. Moreover, if  $\varphi\in\Gamma^\infty$, then $\varphi$ locally coincides with sections belonging $\Gamma_2^\infty$, and therefore $A\varphi$ locally coincides with sections belonging $\Gamma_2^0(\Lambda)$. The latter implies that $A(\Gamma^\infty)\subseteq \Gamma^0(\Lambda)$. The last claim is straightforward.

\end{proof}
\begin{Remark}
{\rm
The previous result depends only on the continuous structure of the field of Hilbert spaces. In other words, the same proof shows an analogous result if $H\to\Lambda$ is a continuous field of Hilbert spaces and $\Gamma^\infty_2$ is replaced by any dense domain in $\mathcal{H}$ invariant by multiplications of functions belonging to $C_c(\Lambda)$. 
}
\end{Remark}

Reduction theory implies the existence of a diagonalization for every self-adjoint operators. More precisely, if $H_0$ is a self-adjoint operator on a Hilbert space $\Hi$, then there is a measurable field of Hilbert spaces $\{H(\lb)\}$ over the spectrum $\sigma(H_0)$ of $H_0$ and a unitary operator $T:\Hi\to\int^\oplus_{\sigma(H_0)}\Hi(\lb)\de\eta(\lb)$ such that $T f(H_0)\varphi(\lb)=f(\lb)T\varphi(\lb)$ for every measurable function $f$ on $\sigma(H_0)$, where $f(H_0)$ denotes the operator defined by the functional calculus and $\eta$ is the so called scalar spectral measure of $H_0$. Moreover, the self-adjoint operators $A$ strongly commuting with $H_0$ (quantum constants of motion) admit a decomposition through $T$. The examples of smooth fields of Hilbert spaces, direct integral and smooth fields of operators that we will consider in section \ref{EX} comes from the diagonalization of certain self-adjoint operator.    

\begin{Corollary}\label{Hor}
Assume $\Lambda\subset \R$ is open and let $H_0=\int_\Lambda^\oplus \lb\de\eta(\lb)$. If $A:\Gamma^\infty_2\to\Gamma^0_2$ is essentially self-adjoint and its closure $\overline{A}$ strongly commutes with $H_0$, then $A\in\AA_2^0$. 
\end{Corollary}
 
Let $A=\{A(\lb)\}$ be a field of operators such that the domain of each $A(\lb)$ contains $\Hi^\infty(\lb)$, $A(\Gamma^\infty_2)\subseteq\Gamma^1(\Lambda)\cap\Hi$ and the same properties hold for $A^*$. Since every section in $\Gamma^1(\Lambda)$ locally coincides with sections belonging $\Gamma^1(\Lambda)\cap\Hi$, it follows that $A$ belongs to $\AA^1$. Therefore  $\hat\nabla_X (A)(\Gamma^\infty)\subseteq\Gamma^0(\Lambda)$, for each $X\in\operatorname{Vect}(\Lambda)$. In order to show that the restriction of $\hat\nabla_X (A)$ to $\Gamma^\infty_2$  defines an operator on $\Hi$ (i.e. $\hat\nabla_X (A)(\Gamma_2^\infty)\subseteq\Gamma^0(\Lambda)\cap\Hi$), we need to make stronger assumptions.      
 
\begin{Definition}\label{smoothdir}
Assume that $\Gamma^\infty_2$ is invariant by $\nabla$. We denote by $\Gamma^n_2(\Lambda)$ the space formed by all the sections $\varphi\in\Gamma^n(\Lambda)$ such that $\nabla_{X_1}\cdots\nabla_{X_k}\varphi\in\Hi$, for every $0\leq k\leq n$ and $X_1,\cdots, X_k\in\operatorname{Vect}(\Lambda)$. We denote by $\AA^n_2$ the space formed by all the fields of operators $A=\{A(\lb)\}$ such that
\begin{enumerate}
\item[a)] the domain of $A(\lb)$ and $A^*(\lb)$ contains $\Hi^\infty(\lb)=\{\varphi(\lb)\mid \varphi\in\Gamma^\infty\}.$
\item[b)]  $A(\Gamma_2^\infty)\subseteq \Gamma_2^n(\Lambda)$ and $A^*(\Gamma^\infty_2)\subseteq \Gamma^n_2(\Lambda)$.
\end{enumerate}
\end{Definition} 

\begin{Remark}
{\rm 
By definition we have the following properties.
\begin{enumerate}
    \item[i)] $\Gamma^n_2(\Lambda)\subseteq\Gamma^n(\Lambda)$, $\Gamma^n_2(\Lambda)\subseteq\Gamma^{n-1}_2(\Lambda)$ and $\nabla_X(\Gamma^n_2(\Lambda))\subseteq\Gamma^{n-1}_2(\Lambda)$, for all $X\in\operatorname{Vect}(\Lambda)$.
    \item[ii)] $\AA^n_2\subseteq\AA^n$, $\AA^n_2\subseteq\AA^{n-1}_2$  and $\hat\nabla_X(\AA^n_2)\subseteq\AA^{n-1}_2$, for all $X\in\operatorname{Vect}(\Lambda)$. Compare the latter inclusions with corollary \ref{cor21}.
\end{enumerate}
}
\end{Remark}

Let us return to the trivializable case and assume $\de\eta$ is a volume form. Let $X$ be a complete vector field on $\Lambda$ and $r_t$ its flow (i.e. $r_t(\lb)$ is the integral curve of $X$ passing through $\lb$ at $t=0$). Also, let $J_t$ be the Jacobian of $r_t$ with respect to $\de\eta$. Therefore, the map $W_t:L^2(\Lambda,V)\to L^2(\Lambda,V)$ given by 
\begin{equation}\label{Ut}
W_t f(\lb)=\sqrt{|J_t|}f(r_{t} \lb),
\end{equation}
for $f\in L^2 (\Lambda,V)$, is a unitary strongly continuous one parameter group. It is easy to show that the infinitesimal generator of $W_t$ is 
$$
H_X=-i(X+\frac{1}{2}\text{div}X)
$$
acting on $L^2(\Lambda,V)$. The following result is a consequence of the latter fact.

\begin{proposition}\label{HamX}
Let $\de\eta$ be a volume form on $\Lambda$ and $H\to\Lambda$ be a trivializable smooth field of Hilbert spaces. If $X$ is a complete vector field on $\Lambda$, then the operator $-i(\nabla_X+\frac{1}{2}\mathrm{div}(X))$ is self-adjoint on $\int_{\Lambda}\Hi(\lb)\de\eta(\lb)$, where the divergence is computed with respect to $\de\eta$.
\end{proposition}

Let $r_t^*:C^n(\Lambda,V)\to C^n(\Lambda,V)$ given by $[r^*_t f](\lb)=f(r_t\lb)$, for each $n\in\mathbb{N}\cup\{\infty\}$. Define $R_t=T^*r^*_tT$. If $A\in\AA^1$, then $\langle R_t A R_{-t}\varphi(\lb),\psi(\lb)\rangle_\lb=\langle \hat A(r_t\lb)T\varphi(\lb),T\psi(\lb)\rangle_V$. Therefore, 
$$
\frac{\de}{\de t}\langle R_t A R_{-t}\varphi(\lb),\psi(\lb)\rangle_\lb|_{t=s}=\langle R_s\hat\nabla_X (A)R_{-s} \varphi(\lb),\psi(\lb)\rangle_\lb
$$
The following result asserts that under suitable conditions the previous pointwise smoothness implies weakly smoothness in the direct integral. The latter turn to be the main tool to explicitly compute $\hat\nabla_X(A)$ in the example of section \ref{EX}.

\begin{Theorem}\label{pointw}
Let $A\in\AA^1$, $X\in\operatorname{Vect}(\Lambda)$ and $r_t$ the one parameter flow of $X$.  Define $R_t:\Gamma^\infty(\Lambda)\to\Gamma^\infty(\Lambda)$ by $R_t=T^*r^*_tT$, where $[r^*_t f](\lb)=f(r_t(\lb))$. Assume that $R_t(\Gamma^\infty)\subseteq\Gamma^\infty$ and let $\varphi,\psi\in\Gamma^\infty$ with compact support such that $\varphi(\lb)\in T^*_\lb (V^\infty)$, for every $\lb\in\Lambda$. Then the map $t\to\langle R_t A R_{-t}\varphi,\psi\rangle$ is differentiable and
$$
\frac{\de}{\de t}\langle R_t A R_{-t}\varphi,\psi\rangle|_{t=s}=\langle R_s\hat\nabla_X (A)R_{-s} \varphi,\psi\rangle
$$
\end{Theorem}
\begin{proof}
It is enough to show that the map $s\to \langle R_s\hat\nabla_X (A)R_{-s} \varphi,\psi\rangle$ is continuous and 
$$
\langle (R_t A R_{-t}-A)\varphi,\psi\rangle=\int_0^t\langle R_s\hat\nabla_X (A)R_{-s} \varphi,\psi\rangle\de s.
$$
Since $A\in\AA^1$, the latter identity holds pointwise, i.e.
$$
\langle (R_t A R_{-t}-A)\varphi(\lb),\psi(\lb)\rangle_\lb=\int_0^t\langle R_s\hat\nabla_X (A)R_{-s} \varphi(\lb),\psi(\lb)\rangle\de s.
$$
For a fix $\lb$, since $r_t$ is a diffeomorphism and $\varphi(\lb)\in T^*_\lb(V^\infty)$, the map 
$$
s\to\langle R_s\hat\nabla_X (A)R_{-s} \varphi(\lb),\psi(\lb)\rangle=\langle \hat\nabla_X (A)(r_s(\lb)) \varphi(\lb),\psi(\lb)\rangle
$$
is continuous. For a fix $s$, since $\hat\nabla_X(A)\in\AA^0$, the section $R_s\hat\nabla_X (A)R_{-s} \varphi$ belongs to $\Gamma^0(\Lambda)$, thus the latter map is also continuous with respect to $\lb$ (so it is continuous and compact supported on $[0,t]\times \Lambda$). In particular, the map  $s\to \langle R_s\hat\nabla_X (A)R_{-s} \varphi,\psi\rangle=\int_\Lambda \langle R_s\hat\nabla_X (A)R_{-s} \varphi(\lb),\psi(\lb)\rangle_\lb\de\eta(\lb)$ is also continuous. Moreover, Fubini's theorem implies that
$$
\langle (R_t A R_{-t}-A)\varphi,\psi\rangle=\int_{\Lambda} \langle (R_t A R_{-t}-A)\varphi(\lb),\psi(\lb)\rangle_\lb\de\eta(\lb)=\int_{\Lambda}\int_0^t\langle R_s\hat\nabla_X (A)R_{-s} \varphi(\lb),\psi(\lb)\rangle\de s\de\eta(\lb) 
$$
$$
=\int_0^t\int_{\Lambda}\langle R_s\hat\nabla_X (A)R_{-s} \varphi(\lb),\psi(\lb)\rangle\de\eta(\lb)\de s=\int_0^t\langle R_s\hat\nabla_X (A)R_{-s} \varphi,\psi\rangle\de s.
$$
\end{proof}

\section{An important example and Canonical Quantization}\label{EX}

 Let us consider the following simple but yet fundamental example: Let $Q^2$ be the multiplication operator on $L^2(\R^n)$ corresponding to the function
$\phi(q)=||q||^2$.  As we previously mentioned, every self-adjoint operator admits a diagonalization. The map  $T:L^2(\R^n)\to\int_{(0,\infty)}^\oplus L^2(\mathbb{S}^{n-1})\de\lb$ given by 
\begin{equation}\label{triv}
T f(\lb,z)=2^{-1/2}\lb^{\frac{n-2}{4}}f(\sqrt{\lb}z)
\end{equation}

is a diagonalization of $Q^2$ (for instance, see lemma 3.6 in \cite{T} for the more general case $\phi(q)=\alpha(\|q\|)$, and for an arbitrary submersion $\phi:\R^n\to \R^k$ see theorem 5.2 in \cite{B}).

Moreover, $T$ can be regarded as an smooth trivialization considering  $\Hi(\lb)=L^2(\mathbb{S}^{n-1}_{\sqrt\lb},\mu_\lb)$, where $\mu_\lb=2^{-1/2}\sqrt{\lb}\eta_\lb$ and $\eta_\lb$ is the canonical measure on $\mathbb{S}^{n-1}_{\sqrt\lb}$. In particular, the restriction of $T$ to $\Hi(\lb)$ defines a unitary operator onto $L^2(\mathbb{S}^{n-1})$. The latter fact allows to identify $f\in C^\infty(\R^n)$ with a section of the field of Hilbert spaces  $\{(0,\infty)\ni\lb\to L^2(\mathbb{S}^{n-1}_{\sqrt\lb},\mu_\lb)\}$  through the restriction $f(\lb)=f|_{\mathbb{S}^{n-1}_{\sqrt\lb}}$. Under that identification, the action of any  $a\in C^\infty(\Lambda)$ on a section $\varphi$ is given by $a\varphi(q)=a(\phi(q))\varphi(q)$.

We will use $T$ to pullback the trivial smooth structure on the trivial field of Hilbert spaces $\hat H=(0,\infty)\times L^2(\mathbb{S}^{n-1})$ into the field of Hilbert spaces $\{(0,\infty)\ni\lb\to L^2(\mathbb{S}^{n-1}_{\sqrt\lb},\mu_\lb)\}$, i.e. we will consider $\nabla_X=T^{-1}X T$. Taking derivatives in equation \eqref{triv}  we obtain the following result.

\begin{proposition}\label{con1}
The map $\nabla_X:C_c^\infty(\R^n)\to C_c^\infty(\R^n)$ given by
$$
\nabla_X(\varphi)=\tilde X(\varphi)+\frac{n-2}{4}\phi^{-1}\tilde X(\phi)\varphi
$$
defines an smooth structure on the field of Hilbert spaces  $\{(0,\infty)\ni\lb\to L^2(\mathbb{S}^{n-1}_{\sqrt\lb},\mu_\lb)\}$, 
where $\tilde X=D\Psi (X)$ and $\Psi:(0,\infty)\times \mathbb{S}^{n-1}\to\R^n\setminus\{0\}$ is the diffeomorphism given by $\Psi(\lb,z)=\sqrt{\lb}z$. Moreover, the map $T$ defined by equation \eqref{triv} is a trivialization of the latter smooth field.
\end{proposition}

\begin{Remark}
{\rm
In the previous proposition one could define the connection  over the sections space $C^\infty(\R^n)$, but we prefer to take $\Gamma^\infty=C_c^\infty(\R^n)$, because the fields of operators that we will consider later are defined on $C_c^\infty(\R^n)$. Moreover, we also have that $\Gamma_2^\infty=C_c^\infty(\R^n)$ and $C^\infty(\R^n)\subset \Gamma^\infty(0,\infty)$. 
}
\end{Remark}
\begin{Remark}\label{flowrtilde}
{\rm
Notice that $\de\phi(\tilde X)=X$, and $\tilde X$ is the only vector field normal to each sphere satisfying that identity. Moreover, if $\tilde r_t$ is the flow of $\tilde X$, then $\phi(\tilde r_t q)=r_t\phi(q)$.
}
\end{Remark}
Let us consider the operator $R_t=T^*r_t^* T$. A direct computation shows that 
$$
R_t \varphi(q)=a(t,q)\varphi(\tilde r_{t}q),
$$
where
$$
a(t,q)=\frac{n-2}{4}\phi(q)^{\frac{2-n}{4}}(r_t\circ\phi(q))^{\frac{n-2}{4}}=\frac{n-2}{4}\phi(q)^{\frac{2-n}{4}}(\phi\circ \tilde r_t(q))^{\frac{n-2}{4}}\,.
$$
Notice that, if $W_t$ is the one-parameter unitary group defined by equation \eqref{Ut}, then $T^*W_t T\varphi(q)=\sqrt{(|J_t|}\circ \phi) a(t,q)\varphi\circ\tilde r_ t(q)$. However, there is a unique one-parameter unitary group of that form and it is given by $\tilde W_t \varphi=\sqrt{|\tilde J_t|}\varphi\circ\tilde r_t$, where $\tilde J_t$ is the Jacobian of $\tilde r_t$. Therefore,
$$
a=\sqrt{|\tilde J_t J_{-t}\circ\phi|}.
$$
Moreover, since $\frac{\de}{\de t}R_t \varphi|_{t=0}=\nabla _X \varphi$, we obtain the following (geometrical) expression for our connection:
\begin{equation}\label{con2}
\nabla_X(\varphi)=\tilde X(\varphi)+\frac{1}{2}(\text{div}(\tilde X)-\text{div}( X)\circ\phi)\varphi.    
\end{equation}
The last expression defines a connection in a much general framework, but we will address that problem in a forthcoming article (for instance, the latter expression defines a connection if $\phi:M\to N$ is a submersion and $\tilde X$ is the unique vector field normal to each $M_\lb=\phi^{-1}(\lb)$ such that $D\phi(\tilde X)=X$, where $M$ and $N$ are Riemannian manifolds).

Note that the operator $-i(\nabla_X+\frac{1}{2}\text{div}(X)\circ\phi)$ considered in proposition \ref{HamX} coincides with $H_{\tilde X}=-i(\tilde X+\frac{1}{2}\text{div}(\tilde X))$ and it is the infinitesimal generator of $\tilde W_t$.

In the next subsection we will compute $\hat\nabla_X(A)$ for suitable $A\in\AA^1$ and to do that we will need the following lemma.
\begin{Lemma}\label{X_0}
Let $X_0$ be the vector field given by $X_0(\lambda)=2\lb\frac{\partial}{\partial\lb}$ . Then $\tilde X_0(q)=\sum q_j\frac{\partial}{\partial q_j}$, 
\begin{equation}\label{R}
  R^0_t \varphi(q)=e^{\frac{n-2}{2}t}\varphi(e^t q)
 \end{equation}
 and 
 $$
 \tilde W^0_t\varphi(q)=e^{\frac{n}{2}t}\varphi(e^t q).
 $$
In particular, $\varphi\in C^\infty(\R^n)$ is an horizontal section if and only if
$$
\varphi(\lb q)=\lb^{-\frac{(n-2)}{2}}\varphi(q),\quad \forall q\in\R^n, q\neq 0,\forall\lb >0.
$$
\end{Lemma}

\begin{proof}
Since $\de\phi(\frac{\nabla\phi}{||\nabla\phi||^2})=1$ and $\de\phi\circ D\Psi=I$,  it follows that $D\Psi(X_0)=\tilde X_0:=\sum q_j\frac{\partial}{\partial q_j}$. Clearly
$r^0_t a(\lb)= a(e^{2t}\lb)$, for every $a\in C^\infty(0,\infty)$. The rest of the proof is a straightforward computation.
 \end{proof}

\subsection{Smooth Fields of Operators Coming from Canonical Quantization}\label{ehfo}

In this subsection we will construct fields of operators over the field of Hilbert spaces $\{(0,\infty)\ni\lb\to L^2(\mathbb{S}^{n-1}_{\sqrt{\lb}},\mu_\lb)\}$ and look for conditions to guarantee that such fields are smooth, or even horizontal. 

Recall that the self-adjoint operators on $L^2(\R^n)$ admitting a decomposition through $T$ are the operators that strongly commutes with $Q^2$ (i.e. the quantum constant of motion of $Q^2$). Such decomposition also holds for bounded (not necessarily self-adjoint) operators strongly commuting with $Q^2$. Moreover, theorem \ref{desc1} also provides conditions to guarantee such decomposition within our smooth framework.

The operators that we shall consider arise from canonical quantization, i.e. they are of the form $\Op(u)$, where $u$ is a ``reasonable function'' on $\R^{2n}$ and $\Op$ is the canonical Weyl calculus \cite{Wey,Fol}. 
Formally, we will prove that Weyl Calculus maps classical constant of motion of $\phi(q)=\|q\|^2$ (where $\phi$ is seen as a function on $\R^{2n}$ independent of the momentum variable) into quantum constant of motion of $Q^2$, so in that way we will obtain an important set of fields of operators where we might apply our connection $\hat\nabla$. 

A classical constant of motion of a classical Hamiltonian $h\in C^\infty(\R^{2n})$ is an smooth function $u\in C^\infty(\R^{2n})$ such that $\{h,u\}=0$, where $\{\cdot,\cdot\}$ is the canonical Poisson bracket on $\R^{2n}$. It is well-known that $u$ is a classical constant of motion if and only if $u\circ\alpha_t=u$, where $\alpha_t$ is the Hamiltonian flow of $h$. It is straightforward to check that, if $h(q,p)=\phi(q)$, then $\alpha_t(q,p)=(q,p+t\nabla\phi(q))=(q,p+2tq)$. Note that in this case $\alpha_t$ is linear. 

Weyl calculus is meant to map real smooth functions on $\R^{2n}$ into self-adjoint operators on $L^2(\R^n)$, but this is not always the case. Indeed, Weyl calculus is a continuous isomorphism $\Op:S'(\R^{2n})\to\mathcal{B}(S(\R^n),S'(\R^n))$, where $S(\R^m)$ is the Schwartz space endowed with its canonical locally convex topology, $S'(\R^m)$ is the topological dual of $S(\R^m)$, i.e the space of tempered distributions and $\mathcal{B}(S(\R^n),S'(\R^n))$ is endowed with the strong operator topology. 

One of the main properties of Weyl calculus is the following identity: if $\varphi,\psi\in S(\R^n)$, then 
$$
\langle\Op(u)\varphi,\psi\rangle=\langle\varphi,\Op(\overline u)\psi\rangle,
$$
where the complex conjugation of distributions is defined by $\overline u(\varphi)=\overline{u(\overline{\varphi}})$ (Proposition 2.6 in \cite{Fol}).

The following property of $\Op$ is our main tool in this section. Let $Sp(2n)$ be the symplectic group, i.e. the group of linear symplectomorphism on $\R^{2n}$; there is a map $m:Sp(2n)\to \mathcal{U}(L^2(\R^n))$ (called the metaplectic representation) such that for any $S\in Sp(2n)$ and every $u\in S'(\R^{2n})$, we have that
\begin{equation}\label{meta}
\Op(u\circ S^*)=m(S)\Op(u)m(S)^{-1}.
\end{equation}

For a detailed presentation of the metaplectic representation see chapter 4 in \cite{Fol}. In particular, see theorem 2.15 for the previous identity. 

Note that equation \eqref{meta} makes sense  because the pull back by $S^*$ maps $S(\R^{2n})$ into itself and it can be extended to an isomorphism from $S'(\R^{2n})$ into itself.

\begin{Definition}
Let $h\in C^\infty(\R^{2n})$ and assume that its Hamiltonian flow $\alpha_t$ is linear and defined for any $t\in\R$. We say that $u\in S'(\R^{2n})$ is a tempered constant of motion if $\alpha_t^*u=u$, for all $t\in\R$. 
\end{Definition}
\begin{Theorem}\label{wpcm}
Let $u$ be a tempered constant of motion of $h(q,p)=\|q\|^2$. 
\begin{enumerate}
\item[a)]  $[\Op(u),e^{itQ^2}]=0$ on $S(\R^n)$, for all $t\in\R$. 
\item[b)] If $\Op(u)$ sends $S(R^n)$ into itself, then $[\Op(u),Q^2] =0$ on $S(R^n)$.
\item [c)] If $\Op(u)$ is bounded or essentially self-adjoint on $S(\R^n)$, then it strongly commutes with $Q^2$.
\item[d)] $\Op(u)\in\AA^0_2$ if and only if $\Op(u)(C^\infty_c(\R^n))\subseteq\Gamma_2^0(0,\infty)$ and $\Op(\overline u)(C^\infty_c(\R^n))\subseteq\Gamma_2^0(0,\infty)$.
\end{enumerate}
\end{Theorem}

\begin{proof}
Equation 4.25 in \cite{Fol} implies that $m(\alpha_t^{*}) = e^{2itQ^2}$. Therefore, equation \eqref{meta} implies that
$$
\Op(u)=e^{itQ^2}\Op(u)e^{-itQ^2}.
$$
and this shows $a)$. When $\Op(u)$ sends $S(R^n)$ into itself, taking strong derivatives in the previous equality we obtain $b)$. Clearly, $a)$ implies the bounded case in $c)$. If $\Op(u)$ is essentially self-adjoint on $S(\R^n)$, let $\mathfrak{D}$ be the domain of $\overline{\Op(u)}$ (the closure of $\Op(u)$). We will prove that $e^{itQ^2}\mathfrak{D}\subseteq \mathfrak{D}$ and $e^{itQ^2}\overline{\Op(u)}=\overline{\Op(u)}e^{itQ^2}$ on $\mathfrak{D}$. Let $f\in \mathfrak{D}$. Then, there is a sequence $f_n\in S(\R^n)$ such that $f_n\to f$ and $\Op(u)f_n\to \overline{\Op(u)}f$. Thus, $e^{itQ^2}f_n\to e^{itQ^2}f$ and 
$$
\Op(u)e^{itQ^2}f_n=e^{itQ^2}\Op(u)f_n\to e^{itQ^2}\overline{\Op(u)}f. 
$$
Since $\overline{\Op(u)}$ is closed and $g_n:=e^{itQ^2}f_n$ is convergent in the graph topology, we have that $e^{itQ^2}f\in \mathfrak{D}$ and 
$$
\overline{\Op(u)}e^{itQ^2}f=e^{itQ^2}\overline{\Op(u)}f.
$$
This implies that $\overline{\Op(u)}$ and $Q^2$ strongly commute.

Let us show d). According to theorem \ref{desc1}, we only need to prove that $a(Q^2)\Op(u)=\Op(u)a(Q^2)$ on $C^\infty_c(\R^n)$, for every $a\in C_c^\infty(0,\infty)$. Let $\hat a$ be the Fourier transform of $a$. For each $\varphi,\psi\in C_c^\infty(\R^n)$, the Fubini's theorem implies that
$$
\int_{-\infty}^\infty \hat a(t)\langle e^{itQ^2}\Op(u)\varphi,\psi\rangle\de t=\int_0^\infty\left(\int_{-\infty}^\infty \hat a(t)e^{it\lb}\de t\right)\langle [\Op(u)\varphi](\lb),\psi(\lb)\rangle_\lb\de\lb=\langle a(Q^2)\Op(u)\varphi,\psi\rangle.
$$
Using a) and repeating the same argument, we obtain that
$$\int_{-\infty}^\infty \hat a(t)\langle e^{itQ^2}\Op(u)\varphi,\psi\rangle\de t=\langle a(Q^2)\varphi,\Op(\overline u)\psi\rangle=\langle \Op(u)a(Q^2)\varphi,\psi\rangle.$$
\end{proof}
\begin{Remark}
{\rm It is well-known that if $u$ belongs to a global H\"ormander class $S^m_{\rho,\delta}(\R^{2n})$, then $\Op(u)(S(\R^n))\subseteq S(\R^n)$ (for instance, see theorem 2.21 in \cite{Fol}). Therefore, if $u$ is a constant of motion of $h(q,p)=\|q\|^2$ belonging to a H\"ormander class $S^m_{\rho,\delta}(\R^{2n})$, then $\Op(u)\in\AA^0_2\cap\AA^{\infty}$. Similarly, if $u$ belongs to a local Hörmander class (see definition 1.1 in \cite{Tay}), it is not difficult to prove that $\Op(u)(C_c^\infty(\mathbb{R}^n))\subseteq S (\mathbb{R}^n)$, for instance repeating the proof of theorem 2.21 in \cite{Fol}. Therefore, if $u$ is a tempered constant of motion of $h(q,p)=\|q\|^2$ belonging to a local H\"ormander class, then $\Op(u)\in\AA^0_2\cap\AA^{\infty}$.
}
\end{Remark}
The main result of this section is a formula to compute $\hat\nabla_X\Op(u)$ (see theorem \ref{Xu} below), but in order to understand it, we need to recall some well-known construction in symplectic geometry. For $\tilde X\in \operatorname{Vect}(\R^n)$, we define $h_{\tilde X}\in C^\infty(\R^{2n})$ by 
\begin{equation}\label{J}
 h_{\tilde X}(q,p)=\langle \tilde X(q),p\rangle,    
 \end{equation}
where $\langle\cdot,\cdot\rangle$ is the duality between the tangent and the cotangent plane at $q$. Equivalently, if $\tilde X=\sum_j^n a_j(q)\frac{\partial}{\partial p_j}$, then $h_{\tilde X}(q,p)=\sum_j^n a_j(q)p_j$.

We will denote by $\hat X$ the Hamiltonian vector field corresponding to $h_{\tilde X}$, i.e. $\hat X(u)=\{h_{\tilde X},u\}$. Recall that the flow of $\tilde X$ is denoted by $\tilde r_{t}$ (see remark \ref{flowrtilde}). It is easy to show the flow $\hat r_t$ of $\hat X$ is the canonical lift of $\tilde r_t$, i.e. $\hat r_t(q,p)=(\tilde r_t(q),\tilde r^*_{-t}(p))$.

It is also well-known that $\Op(h_{\tilde X})=H_{\tilde X}=-i(\tilde X+\frac{1}{2}\text{div}(\tilde X))$. In particular, if $\tilde X=D\Psi(X)$, then $\Op(h_{\tilde X})=-i(\hat\nabla_X+\frac{1}{2}\text{div}(X))$.

When $X_0(\lambda)=2\lb\frac{\de}{\de\lb}$, we have that $h_0(q,p):=h_{\tilde X_0}(q,p)=\sum q_j p_j$ and $\hat r_t^{0}(q,p)=(e^t q,e^{-t}p)$. In particular, $\hat r_t^{0}$ is linear. Moreover, equation 4.24 in \cite{Fol} implies that
$$m(\hat r_t^{0})=\tilde W_t^{0}.$$ 
The latter identity allow us to compute $\tilde W_t^{0}\Op(u)\tilde  W_{-t}^{0}$ using the metaplectic representation, and it is the main reason why we are considering the vector field $X_0$. Furthermore, since  $e^{t}R_t^{0}=\tilde W_t^{0}$ (lemma \ref{X_0}), we obtain the identity 
\medskip
$$
R_t^{0} \Op(u) R_{-t}^{0}=\Op( (\hat r_t^{0})^* u).
$$
Using a Taylor expansion, it is straightforward to show that the limit $\lim_{t\to 0}\frac{1}{t} (u\circ \hat r_t^{0}-u)=\{h_0,u\}$ holds in $S'(\R^{2n})$. Since $\Op: S'(\R^{2n})\to B(S(\R^n),S'(\R^n))$ is continuous, theorem \ref{pointw} implies the following remarkable result.
\begin{Theorem}\label{Xu}
 If $u$ is a tempered constant of motion of $h(q,p)=\|q\|^2$, $X=aX_0$  and $\Op(u)\in\AA^1$, then 
\begin{equation}\label{MF}
  \hat\nabla_{X}\Op(u)=a(Q^2)\Op(\hat X_0(u)). 
 \end{equation}
 In particular, if $\Op(u)$ and $\Op(\hat X_0(u))$ belong to $\AA^0_2$, then $\Op(u)\in\AA^1_2$.
\end{Theorem}
\begin{Remark}
{\rm
Since $\{h,h_0\}=2h$, the Jacobi identity implies that $\{h_0,u\}=\hat X_0(u)$ is a tempered constant of motion of $h$.
}
\end{Remark}
\begin{Corollary}\label{cor310}
If $u$ is a tempered constant of motion of $h$, then $\Op(u)$ is horizontal if and only if $u$ is also a tempered constant of motion of $h_0$ and $\Op(u)\in\AA^0_2$.
\end{Corollary}
\begin{proof}
If $\Op(u)$ is horizontal, theorem \ref{Xu} implies that $\Op(\hat X_0(u))=0$. Since $\Op$ is faithful, $u$ is a constant of motion of $h_0$. The converse is a direct consequence of lemma 5.1.1. in \cite{LS}. 
\end{proof}
\begin{Corollary}\label{Cor310}
If $u$ is a tempered constant of motion of $h(q,p)=\|q\|^2$ belonging to a local H\"ormander class $S^m_{\rho,\delta}(\R^{2n})$, then $\Op(u)$ belongs to $\AA^\infty_2$.
\end{Corollary}
\begin{proof}
We already note that $\Op(u)$ belongs to $\AA^0_2\cap\AA^{\infty}$. Since $\{h_0,u\}\in S^{m'}_{\rho,\delta}(\R^{2n})$ with $m'=m+\max(1-\rho,\delta)$, theorems \ref{wpcm} and \ref{Xu} implies our result. 
\end{proof}
The following result is analogue to theorem \ref{Xu}, but taking $h(q,p)=\|p\|^2$. The corresponding operator is the Laplacian $H=-\Delta$, which is fundamental for physical applications, and its diagonalization $\tilde T$ is obtained after composing $T$ defined by equation \eqref{triv} with the Fourier transform $\mathcal F$. In order to make $\tilde T$ a trivialization, we define the connection by $\nabla_X=\tilde T^* X\tilde T$. Also notice that the metaplectic representation $m$ maps the symplectic matrix $\mathfrak J$ into the Fourier transform $\mathcal F$. Moreover, $h\circ \mathfrak J=\tilde h$ and $h_0\circ \mathfrak J=- h_0$, where $\tilde h(q,p)=\|q\|^2$. Since $\mathfrak J$ is a symplectomorphism, $u$ is a constant of motion of $h$ if and only if $u\circ\mathfrak J$ is a constant of motion of $\tilde h$.
\begin{Corollary}\label{Lapl}
 If $u$ is a tempered constant of motion of $h(q,p)=\|p\|^2$, $X=aX_0$  and $\Op(u)\in\AA^1$, then 
 $$
 \hat\nabla_{X}\Op(u)=a(-\Delta)\Op(\hat X_0(u)).
 $$
 In particular, if $\Op(u)$ and $\Op(\hat X_0(u))$ belong to $\AA^0_2$, then $\Op(u)\in\AA^1_2$.
\end{Corollary}
\begin{Remark}
{\rm Through the rest of this article we will only consider fields of operators coming from Weyl quantization. However, there are other ways to construct operators admiting a decomposition through $\tilde{T}$. For example, if $S$ is the scattering operator corresponding to a suitable Schrödinger operator, then $S$ is unitary and it strongly commutes with $-\Delta$, therefore it can be decomposed through $\tilde{T}$ (for instance, see \cite{Yaf}). For the moment, we do not know if our notion of smooth fields of operators and our results can be applied in scattering theory, and we look forward to study this problem in the future. 
}
\end{Remark}
Let us return to the case $H=Q^{2}$. When $X=X_0$, formula \eqref{MF} becomes
\begin{equation}\label{X0}
   \hat\nabla_{X_0}\Op(u)=\Op(\hat X_0(u)).   
\end{equation}

We do not know if the latter equation holds true when we replace $X_0$ by and arbitrary vector field $X$. Such identity is equivalent to
$$
a(Q^2)\Op(\hat X_0(u))=\Op((a\circ\phi) \hat X_0(u)),
$$
where $\phi(q)=\|q\|^{2}$. It is well-known that $\Op$ is not a multiplicative homomorphism. However, because of the particular type of symbols that we are considering, the latter identity might still hold true. Notice that for an arbitrary vector field $X$, the Hamiltonian flow $\hat r_t$ is not necessarily linear, so we cannot use the metaplectic representation. Nevertheless, in order to estimate the difference between the right and left hand of equation \eqref{X0} for an arbitrary $X$, we might use semi-classical theory, introducing Planck's constant dependence and applying Egorov's theorem (for instance, see theorem 11.1 in \cite{Z}). We expect to address the latter problem in the future. 

The following observation relates our results and the problem discussed in the previous paragraph with deformation quantization. Let $\A$ be the Poisson algebra of constants of motion of $h(p,q)=\|q\|^2$. Then, the map $\tilde\nabla:\operatorname{Vect}((0,\infty))\times\A\mapsto\A$, given by $\tilde\nabla_X(u)=\hat X(u)$ satisfies the following properties: 
\begin{itemize}
    \item[a)] $\tilde\nabla_X(\{u,v\})=\{\tilde\nabla_X u,v\}+\{u,\tilde\nabla_{\overline{X}}v\}$.
    \item[b)] $\tilde\nabla_X\tilde\nabla_Y-\tilde\nabla_Y\tilde\nabla_X=\tilde\nabla_{[X,Y]}$.
\end{itemize}
In other words, $\tilde\nabla$ defines a sort of abelian Poisson connection.

We can generalize the previous construction as follows. Let $M$ be a Riemannian manifold and $\phi\in C^\infty(M)$ regular. For each $X\in\operatorname{Vect}(\phi(M))$, let $\tilde X$ be the unique vector field normal to each $M_\lb:=\phi^{-1}(\lb)$ such that $D\phi(\tilde X)=X$, and let $\hat X$ be the Hamiltonian vector field of $h_{\tilde X}$ on $T^*M$. If $\A$ is the Poisson algebra of constants of motion of $h(q,p)=\phi(q)$, then the map 
$\tilde\nabla:\operatorname{Vect}((\phi(M))\times\A\mapsto\A$, given by $\tilde\nabla_X(u)=\hat X(u)$ is well-defined and satisfies a) and b) above as well.

Let $\alpha_t$ be the Hamiltonian flow of $h(q,p)=\phi(q)$. For each $\lb\in\phi(M)$, let $\hat\Sigma_\lb=h^{-1}(\lb)$ be the constant energy submanifold (do not confuse it with $M_\lb$ above) and denote by $\Sigma_\lb$ the orbit space $\hat\Sigma_\lb/\alpha$ endowed with the symplectic structure obtained after applying Marsden-Weinstein-Meyer reduction \cite{MW,Me} (one of the authors showed in \cite{B} that $\Sigma_\lb=T^*M_\lb$). For each $u\in\A$, we can define $u_\lb\in C^\infty(\Si_\lb)$ by $u_\lb([\sigma])=u(\sigma)$, where $\sigma\in\hat\Sigma_\lb$ and $[\sigma]$ denotes the orbit $\sigma$. Then, every $u\in \A$ can be regarded as a section of a field of Poisson algebras over $\phi(M)$ with fibers $C^\infty(\Si_\lb)$. The relation between symplectic connections and deformation quantization has been successfully studied, specially since Fedosov's  paper \cite{Fed}. Fedosov showed that, if $\Sigma$ is a symplectic manifold then the Poisson algebra $C^{\infty}(\Sigma)$ admits a star product. In order to do so, Fedosov used an abelian symplectic connection to glue the canonical Moyal product on each tangent plane. In our framework, if a star product $\star^\lb_\hbar$ is given on each $C^\infty(\Sigma_\lb)$, we would like to adapt Fedosov's ideas to use our Poisson connection to glue those star products into a single one defined on $\A$. If the latter construction works, then we might wonder if there is a relation between such star product and the canonical Moyal product restricted to $\A$. An analogous question was formulated within the framework of canonical quantization and Wigner transforms in \cite{B}. Both questions can be interpreted as problems analogue to the commutation of reduction and quantization problem in geometric quantization theory \cite{GS, Woo}. We shall leave the latter problem open as well.

\subsection{Functions of Angular Momenta as Horizontal Constants of Motion.}\label{AM}

Let us construct some examples of constants of motion. For future references, we will assume for a while that $\phi$ is any smooth function on an smooth manifold $M$. 
 
Let $G$ be a Lie group acting on $M$. Such action induces the Lie algebra homomorphism $\zeta:\mathfrak{g}\to\operatorname{Vect}(M)$ given by
$$
\zeta (X)(q)=\frac{\de}{\de t}\left(\exp(t X)\cdot q\right)|_{t=0},
$$
where $\mathfrak g$ is the Lie algebra corresponding to $G$ and $X\in\mathfrak g$. Also let the $\mathcal J:T^*M\to\mathfrak{g}^*$ the induced moment map given by
$$
\mathcal J[(q,p)](X)=\langle\zeta(X)(q),p\rangle=h_{\zeta(X)}(q,p),
$$
where $\langle\cdot\,,\,\cdot\rangle$ implements the duality between $T_q M$ and $T^*_q M$.

If we endow $\mathfrak{g}^*$ with the coadjoint action of $G$ and the negative of the canonical Lie-Poisson structure, and we lift the action of $G$ to $T^*M$, then $\mathcal J$ turn to be an equivariant Poisson map.

\begin{proposition}\label{LP}
Let $G$ be a Lie group acting on $M$ and $\mathcal J$ be the moment map defined above.
\begin{enumerate}
    \item[i)] If $\phi(g\cdot q)=\phi(q)$, for each  $q\in M$, then $a\circ \mathcal{J}$ is constant of motion of  $h(q,p)=\phi(q)$, for any $a\in C^\infty(\mathfrak{g}^*)$.
    \item[ii)] If $\tilde Y$ is a vector field on $M$ such that $[\tilde Y,\zeta(X)]=0$ for every $X\in\mathfrak{g}$, then $a\circ \mathcal{J}$ is constant of motion of $h_{\tilde Y}$, for any $a\in C^\infty(\mathfrak{g}^*)$. 
\end{enumerate}

\end{proposition}

\begin{proof}
Suppose that $\phi(g\cdot q)=\phi(q)$, for each $j$ and each $q\in M$. Thus, $\zeta(X)$ is tangent to each $M_\lb$. Moreover, for each $t\in\R$ and $(q,p)\in T^*M$, since $\nabla\phi$ is normal to each $M_\lb$, we have that
$$
[\mathcal J\circ\alpha_t(q,p)](X)=[\mathcal J(q,t\de\phi(q)+p)](X)=\langle\zeta(X)(q),t\de\phi(q)+p\rangle=
$$
$$
\langle\zeta(X)(q),p\rangle=\mathcal J(q,p)(X).
$$
In other words $\mathcal J\circ\alpha_t=\mathcal J$ and this implies i). For the second part, we shall prove that $\mathcal J\circ\hat r_t=\mathcal J$, where $\hat r_t$ is the Hamiltonian flow of $h_{\tilde Y}$. Notice that 
$[\mathcal J\circ\hat r_t(q,p)](X)=(h_{\zeta(X)}\circ\hat r_t)(q,p)$, therefore $\mathcal J\circ\hat r_t=\mathcal J$ if and only if each $h_{\zeta(X)}$ is a constant of motion of $h_{\tilde Y}$. Our result follows from the identity $\{h_{\tilde Y},h_{\zeta(X)}\}=h_{[\tilde Y,\zeta(X)]}$ (equation II 3.11 in \cite{La}).
\end{proof}

Let us return to our example $\phi(x)=\|x\|^2$. Also, let $G=O(n)$ be the orthogonal group acting canonically on $\R^n$. Then clearly conditions i) in the previous proposition is satisfied. Since $\tilde r^0_t(q)=e^{t}q$ and rotations commutes  with dilations, condition ii) is satisfied for $\tilde Y=\tilde X$, for any $X\in\operatorname{Vect}(0,\infty)$.
\begin{Corollary}\label{amh}
Let $O(n)$ act on $\R^n$ canonically and $\mathcal J:\R^{2n}\to\mathfrak{so}(n)^*$ the corresponding moment map. Then $a\circ \mathcal J$ is a classical constant of motion of $h(q,p)=\|q\|^2$ and of $h_{\tilde X}$, for any $X\in\operatorname{Vect}(0,\infty)$ and $a\in C^\infty(\mathfrak{so}(n)^*)$. In particular, if $a\circ\mathcal J$ is a tempered distribution belonging to some local H\"ormander class $S^m_{\rho,\delta}(\R^{2n})$, then $\Op(a\circ\mathcal J)$ is a horizontal field of operators.
\end{Corollary}

\begin{Remark}
{\rm
Let $l_{i,j}(q,p)=q_i p_j- q_jp_i$ and $L_{ij}=\Op(l_{ij})=q_i\frac{\partial}{\partial q_j}-q_j\frac{\partial}{\partial q_i}$. Then $l_{ij}=h_{\zeta(X_{ij})}$, where $X_{ij}$ is the element of the canonical basis of $\mathfrak{so}(n)$ corresponding to the infinitesimal generator of the clockwise rotation on the plane in $\R^n$ generated by $e_i$ and $e_j$. The functions $l_{i,j}$ and the operator $L_{ij}$ are called the classic and quantum angular momenta coordinates respectively.  Notice that defining a polynomial of the family of operators $L_{ij}$ is not at all trivial, because the operators $L_{ij}$ do not commute (they come from a representation of $\mathfrak{so}(n)$, so they satisfy the same commutation relations than the corresponding vectors $X_{ij}$). The latter is an angular momenta version of the canonical ordering problem for the operators position and momentum. Weyl calculus can be interpreted as a symmetric solution of the canonical ordering problem, and our result suggest that it is also a convenient solution of the angular momenta ordering problem. 
}
\end{Remark}
\begin{Remark}
{\rm
The spectral analysis of the operators $L_{ij}$ and $L^2=\sum_{i<j}L^2_{ij}$ (the total angular momentum operator) are usually presented in any course on quantum mechanics. Using polar coordinates, one can show that the spectra of those operators coincide with the spectra of their restrictions to $L^2(\mathbb{S}^{n-1})$. Our previous result implies that the same conclusion holds for a much larger class of operators.
}
\end{Remark}

\section{Locally uniformly bounded smooth field of operators and smooth fields of $C^*$-algebras.}\label{LB}

In this subsection, we will study how the smoothness of a field of operators $A\in\AA^n$ interacts with the continuity of $A$ as an operator on $\Gamma^n(\Lambda)$. The discussion will lead us to introduce a notion of smoothness for fields of $C^*$-algebras. 

It will be useful to consider the Hilbert $C_0(\Lambda)$-module $\Gamma^0_0(\Lambda)$. Recall that an operator $A:\Gamma^0_0(\Lambda)\to \Gamma^0_0(\Lambda)$ is called adjointable if and only if there is another operator $A^*:\Gamma^0_0(\Lambda)\to \Gamma^0_0(\Lambda)$ such that $h(A\varphi,\psi)=h(\varphi,A^*\psi)$, for any $\varphi,\psi\in \Gamma^0_0(\Lambda)$. It is well-known that, in contrast with operators on Hilbert spaces, the continuity of $A$ does not guarantee that $A$ is adjointable. On the other hand, if $A$ is adjointable then $A$ is continuous (for instance, see lemma 2.18 in \cite{RW}). Moreover, the space of all adjointable operators form a $C^{*}$-algebra (proposition 2.21 \cite{RW}).

Another important property is that every adjointable operator is given by a field of operators. Indeed, if $A$ is an adjointable operator, then clearly $A(K_0(\lb))\subseteq K_0(\lb)$. Therefore, by lemma \ref{sic} there is an operator $A(\lb)$ on $\Hi(\lb)$ such that $A\varphi(\lb)=A(\lb)\varphi(\lb)$, for all $\varphi\in\Gamma^0_0(\Lambda)$ and $\lambda\in\Lambda$.   

Let us introduce one of the spaces of fields of operators that we will consider in this subsection.  
\begin{Definition}
For each $n\in\mathbb N\cup\{\infty\}$, we denote by $\AA^n_{c}$ the space formed by the field of operators $A\in\AA^n$ such that $A:\Gamma^\infty\to\Gamma^n(\Lambda)$ is continuous with respect to the Fréchet topology of $\Gamma^n(\Lambda)$. We denote by $\AA^0_0$ the space of adjointable operators on $\Gamma^0_0(\Lambda)$.
\end{Definition}

In particular, if $A$ belongs to $\AA^n_c$, then $A$ can be extended to $\Gamma^n(\Lambda)$; we will also denote by $A$ such extension. Let us apply our previous discussion to characterize $\AA^0_c$.

\begin{proposition}\label{zero}
Let $A:\Gamma^\infty\to\Gamma^0(\Lambda)$ be a linear operator. The following statements are equivalent
\begin{enumerate}
    \item[a)] $A\in\AA^0_{c}$
    \item[b)] $A\in\AA^0$ and it is a locally uniformly bounded field of operators.
    \item[c)] $A$ extends to $\Gamma^0(\Lambda)$ and its extension is adjointable.
\end{enumerate}
\end{proposition}
\begin{Remark}
{\rm
We are abusing of the notation in part c) of the previous statement, because we have defined adjointability on $\Gamma^0_0(\Lambda)$, but the corresponding definition is the same on $\Gamma^0(\Lambda)$. 
}

\end{Remark}
\begin{proof}
Let us show that \textit{a)} implies \textit{b)}. Since $A$ is a field of operators, $A|_{\Gamma^0(C)}$ is a well-defined and continuous operator on $\Gamma^0(C)$, for every compact set $C\subseteq\Lambda$.
Moreover, $A K_0(\lb)\subseteq K_0(\lb)$ and lemma \ref{sic} implies  that
$$
\|A(\lb)\varphi(\lb)\|_\lb=\inf\{\|A\varphi-\psi\|\mid\psi\in K_0(\lb)\}\leq \inf\{\|A(\varphi-\psi)\|\mid\psi\in K_0(\lb)\}\leq \|A\|\|\varphi(\lb)\|_\lb,
$$
for all $\varphi\in\Gamma^0(C)$. Since $\Hi(\lb)=\{\varphi(\lb)\mid \varphi\in \Gamma^0(C)\}$, the latter inequality implies \textit{b)}.

Clearly \textit{b)} implies \textit{c)} and \textit{b)} implies \textit{a)}. In order to show that \textit{c)} implies \textit{a)}, notice that since $A$ is adjointable, $\varphi|C=0$ implies that $A\varphi|_C=0$, for every compact $C\subseteq \Lambda$. Therefore the map $A|_C:\Gamma^0(C)\to\Gamma^0(C)$ given by  $A|_C(\varphi|_C)=A(\varphi)|_C$ is well-defined and also adjointable. Since $\Gamma^0(C)$ is a Hilbert module, $A|_C$ is a bounded field of operators and clearly this implies \textit{a)}.
\end{proof}
\begin{Remark}
{\rm
For each $C\subseteq\Lambda$ compact, the previous proof also shows that \begin{equation}\label{knorm}
    \|A\|_C:=\|A|_{\Gamma^0(C)}\|=\sup_{\lambda\in C}\{\|A(\lambda)\|\}
\end{equation} 
The space $\AA^0_{c}$ with the family of seminorms $\|A\|_C$ forms a $C^*$-{\it locally algebra}, i.e. $\AA^0_{c}$ corresponds to a $*$-algebra equipped with a locally convex topology which is Hausdorff, complete and generated by a family of $C^*$-seminorms \cite{Ino}.

}
\end{Remark}

\begin{Remark}
{\rm
Let $A\in\AA^0$ and $C\subseteq\Lambda$ compact. Even if $A|_C$ is not an uniformly bounded field of operators, it is a densely defined operator on $\Gamma^0(C)$. Moreover, it is easy to show that $A|_C$ is closable on $\Gamma^0(C)$. Therefore $A|_C$ is a semiregular operator on $\Gamma^0(C)$ as defined in \cite{AR} (i.e. a closable densely defined operator on a Hilbert $C^{*}$-Module with an adjoint also densely defined). We would like to obtain conditions to guarantee that $A|_C$ is also regular and to study how regularity is related with selfadjointness of such operators but acting on a direct integral, as explained in subsection \ref{IMS}. However, we will not consider the latter problems in this article. 
}
\end{Remark}

Let us look for conditions to guarantee that a given field of operators $A$ belongs $\AA^n_{c}$. Let $\varphi_j, \varphi\in\Gamma^\infty$ such that $\varphi_j\to\varphi$ in $\Gamma^n(\Lambda)$ . Then $A\varphi_j\to A\varphi$ in $\Gamma^n(\Lambda)$ if and only if
\begin{enumerate}
\item[a)] $A\varphi_j\to A\varphi$ in $\Gamma^{n-1}(\Lambda)$
\item[b)] $\nabla_X A\varphi_j=A\nabla_X\varphi_j+\hat\nabla_X(A)\varphi_j\to \nabla_X A\varphi$ in $\Gamma^{n-1}(\Lambda)$, for every $X\in\operatorname{Vect}(\Lambda)$.
\end{enumerate}

Since $\varphi_j\to\varphi$ in $\Gamma^{n-1}(\Lambda)$ and $\nabla_X\varphi_j\to\nabla_X\varphi$ in $\Gamma^{n-1}(\Lambda)$, for every $X\in\operatorname{Vect}(\Lambda)$, we obtain the following result.  

\begin{proposition}\label{inter}
Let $A\in\AA^n$. 
\begin{enumerate}
\item[i)] Assume that $A\in\AA^{n-1}_{c}$. Then $A\in\AA^{n}_{c}$ if and only if $\hat\nabla_X(A):\Gamma^\infty\to\Gamma^{n-1}(\Lambda)$ is continuous when $\Gamma^\infty$ is endowed with the $\Gamma^n(\Lambda)$-topology.
\item[ii)] If $A\in\AA^{n-1}_{c}$ and $\hat\nabla_X(A)\in\AA^{n-1}_{c}$, then $A\in\AA^{n}_{c}$.

\end{enumerate}
In particular, if $\hat\nabla_{X_1}\cdots \hat\nabla_{X_k}A\in\AA^0_{c}$, for every $X_1,\cdots X_k\in\operatorname{Vect}(\Lambda)$ and $0\leq k\leq n$, then $A\in\AA^{n}_{c}$.
\end{proposition}

Notice that, $\AA^n_c\nsubseteq\AA^{n-1}_c$ and $\hat\nabla_X(\AA^n_c)\nsubseteq\AA^{n-1}_c$, for any $X\in\operatorname{Vect}(\Lambda)$. The previous result allow us to introduce the following subspace of $\AA^n_c$, where the corresponding inclusions become true. 
\begin{Definition}\label{def47}
For each $n\in\mathbb N\cup\{\infty\}$, we denote by $\AA^n_{lb}$ the space formed by the fields of operators $A\in\AA^n$ such that the field of operators $\hat\nabla_{X_1}\cdots \hat\nabla_{X_k}A$ is locally uniformly bounded, for any $X_1,\cdots X_k\in\operatorname{Vect}(\Lambda)$ and $0\leq k\leq n$.
\end{Definition}

We could use statements \textit{a)} or \textit{c)} in proposition \ref{zero} to define $\AA^n_{lb}$ , instead of \textit{b)}. 

We can endow $\AA^n_{lb}$ with a family of seminorms, just as it was done for $\Gamma^n(\Lambda)$. More precisely, we define 
$$
||A||_{C,X_1,\cdots X_m}=\sup\{||\hat\nabla_{X_1}\cdots\hat\nabla_{X_m}A(\lb)||:\lb\in C\},    
$$
where $C\subseteq\Lambda$ is compact, $X_1,\cdots,X_m\in\operatorname{Vect}(\Lambda)$ and $m\leq n$. Clearly, $\AA^n_{lb}$ becomes a Fréchet space with the latter family of seminorms.
\begin{Remark}
{\rm
We can also put the natural direct limit topology on the space $\AA^n_{lb}$ (or $\Gamma^n(\Lambda)$), but we will not consider it in this article.
}
\end{Remark}

\begin{Corollary}\label{cor21}
For any $n\geq 1$ and $X\in\operatorname{Vect}(\Lambda)$, $\AA^n_{lb}\subseteq\AA^{n-1}_{lb}$, $\hat\nabla_X(\AA^n_{lb})\subseteq\AA^{n-1}_{lb}$ and the map $\hat\nabla_X: \mathfrak{A}^n_{lb}\to\mathfrak{A}_{lb}^{n-1}$ is continuous.
\end{Corollary}

Also notice that if $A, B\in\AA^n_{lb}$, then $AB\in\AA^n_{lb}$ and the following identity holds
\begin{equation}\label{Lei}
 \hat\nabla_X(AB)=\hat\nabla_X(A)B+A\hat\nabla_X(B). \end{equation}

In the full projectively trivializable case, we obtain stronger properties of the map $\hat T$, defined in theorem \ref{triviau}, if we restrict it to $\AA^n_{lb}$. For instance, since each operator $\hat A(\lb)=T_\lb A(\lb)T^*_\lb$ is bounded and $T^*f\in \Gamma^\infty(\Lambda)\subset \Gamma^n(\Lambda)$,  for any $A\in\AA^n_{lb}$ and $f\in V$, we no longer need to consider the common domain $V^\infty$ defined in \ref{vinf}. Recall that $B(V)_{\ast-s}$ denotes the space of bounded operators on $V$ endowed with the $\ast$-strong topology. We also denote by $C^n_{lb}(\Lambda,B(V)_{\ast-s})$ the space of $n$-times differentiable functions from $\Lambda$ to $B(V)_{\ast-s}$ such that any derivative of order lower or equal to $n$ is continuous and locally bounded.

\begin{proposition}\label{pro410}
Let $T:H\to V$ be a full projective trivialization of the smooth field of Hilbert spaces $H\to\Lambda$. The map $\hat T$, defined by $\hat T A(\lb)=T_\lb A(\lb)T^*_\lb$, is a local isometry of $\AA^n_{lb}$ onto $C^n_{lb}(\Lambda,B(V)_{\ast-s})$. 
\end{proposition}

\begin{proof}
The same proof of theorem \ref{triviau} and the previous comment show that $\hat T$ is a well-defined local isometry. It is enough to show that if $\hat A\in C^n_{lb}(\Lambda,B(V)_{\ast-s})$, then $A=T^*\hat AT\in \AA^n_{lb}$. Equivalently, it is enough to prove that $\hat A(C^\infty(\Lambda,V))\subseteq C^n(\Lambda,V)$. The case $n=0$, follows from the inequality
$$
\|\hat A\varphi(\lb)-\hat A\varphi(\lb_0)\|\leq \|[\hat A(\lb)-\hat A(\lb_0)]\varphi(\lb_0)\|+\|\hat A(\lb)\|\|\varphi(\lb)-\varphi(\lb_0)\|.
$$
Similarly, we can use the canonical proof of Leibniz's multiplication formula to show that $\hat A\varphi\in C^n(\Lambda,V)$, for $n>0$ and every $\varphi\in C^\infty(\Lambda,V)$ (or lemma 5.1.1 in \cite{LS}). 
\end{proof}

\subsection{Smooth fields of $C^*$-algebras}\label{SFCA}

Recall that one of our aims is to propose a definition of a smooth field of $C^*$-algebras. Let us recall the definition of a (upper semi-)continuous field of $C^*$-algebras.

\begin{Definition}\label{contical}
Let $\Lambda$ be a locally compact Hausdorff space and $p: \A \rightarrow \Lambda$ be a field of $C^*$-algebras (i.e $p$ is a surjection such that $\AA(\lb):=p^{-1}(\lb)$ is a $C^{*}$-algebra, for any $\lb\in\Lambda$). A (upper semi-)continuous structure on $p: \A \rightarrow \Lambda$ is given by specifying a $\ast$-algebra of sections $\AA$, closed under multiplication by elements of $C(\Lambda)$, and such that:
 \begin{enumerate}
     \item[i)]  the map $\lb \mapsto\|A(\lb)\|$ is (upper semi-)continuous, for all $A \in \AA$, and
     \item[ii)] $\{A(\lb)\mid A \in \AA\}$ is dense in $\AA(\lb)$, for all $\lb \in \Lambda$.
 \end{enumerate}
\end{Definition}

There are at least two other ways to characterize (upper semi-)continuous fields of $C^*$- algebras and each of them has their own advantages. For instance, the previous definition allows to describe the continuous structure in terms of a space of sections, and this is the perspective that we are following in this article (recall the definition  \ref{def} of smooth fields of Hilbert spaces and compare with the definition \ref{DefC} of smooth field of $C^*$-algebras that we will propose latter). 

As we mentioned in the introduction, the previous definition of (upper semi-)continuous field of $C^*$-algebras is equivalent to the notion of (upper semi-)continuous $C^*$-bundle (see definition C.16 and theorem C.25 in \cite{Wi}). Note that, to obtain a (upper semi-)continuous $C^*$-bundle from $\AA$, it is not necessary to require that $\AA$ is closed under multiplication by elements of $C(\Lambda)$. However, under such condition it is easy to show that $\AA$ is dense in the space of continuous sections $\AA(\Lambda)$ of the $C^{*}$-bundle $\mathcal{A}\to \Lambda$ endowed with the locally uniform convergence topology (see proposition C.24 in \cite{Wi}).

Another equivalent way to describe (upper semi-)continuous field of $C^*$-algebras is through the notion of $C_0(\Lambda)$-algebra (see definition C.1 \cite{Wi}). If $\A\to\Lambda$ is an upper semi-continuous $C^*$-bundle, then the subspace $\AA_0(\Lambda)\subset\AA(\Lambda)$, consisting of continuous sections vanishing at infinity, is a $C_0(\Lambda)$-algebra.
The converse is more interesting and involved (see the proof in proposition C.10 in \cite{Wi}).

\medskip

If $A\in\AA^0_{lb}$, by definition the maps $\lb\mapsto \|A(\lb)\varphi(\lb)\|$ and $\lb\mapsto \|A^*(\lb)\varphi(\lb)\|$ are continuous, for every $\varphi\in\Gamma^0(\Lambda)$. However, the map $\lb\mapsto \|A(\lb)\|$ is not necessarily (upper semi-)continuous. Therefore, $\AA^0_{lb}$ does not define an (upper semi-)continuous structure on the field of $C^*$-algebras $\AA(\lb)=\{A(\lb)\mid A\in\AA^0_{lb}\}$. Notice that $\AA^0_0$ is a $C_0(\Lambda)$-algebra, thus $\AA^0_{lb}$ defines an upper semi-continuous structure but on the field of $C^*$-algebras $\tilde\AA(\lb)=\AA^0_0\big/ I_\lb\AA^0_0$ (and generically this algebra does not coincide with $\AA(\lb)$). 

There are subalgebras $\AA$ of $\AA^0_{lb}$ such that the map $\lb\mapsto \|A(\lb)\|$ is continuous, for any $A\in \AA$. Indeed, we will show that the space of compact operators $\mathbb{K}_0(\Lambda):=\mathbb{K}(\Gamma_0^0(\Lambda))$ on the Hilbert module $\Gamma_0^0(\Lambda)$ satisfies such property. Let us recall the construction of the $C^*$-algebra $\mathbb{K}_0(\Lambda)$. For each $\varphi,\psi\in\Gamma_0^0(\Lambda)$, we define the adjointable operator $|\varphi\rangle\langle\psi|$ by
$$
|\varphi\rangle \langle\psi|(\phi)=\langle \psi,\phi\rangle\varphi
$$
It is straightforward to show that $|\varphi\rangle\langle\psi|^*=|\psi\rangle\langle\varphi|$ and $\||\varphi\rangle\langle\psi|\|=\|\varphi(\lb)\|\|\psi(\lb)\|$. The space of compact operators $\mathbb{K}_0(\Lambda)$ is the closure span of the set of operators of the form $|\varphi\rangle \langle\psi|$\,, with $\varphi,\,\psi\in\Gamma_0^0(\Lambda)$.

The following proposition summarizes the discussion in subsection 10.7.2 \cite{D2}, but it is expressed within the framework of Hilbert modules. For completeness of our presentation we give some details in the proof.

\begin{proposition}\label{exako}
Let $H\to\Lambda$ an smooth field of Hilbert spaces and let $p:\mathcal{K}\to \Lambda$ be the field of $C^*$-algebras with fibers $\mathbb{K}(\Hi(\lb))$, where $\mathbb{K}(\Hi(\lb))$ is the $C^*$-algebra of compact operators on the Hilbert space $\Hi(\lb)$.
Identifying $A\in \mathbb{K}_0(\Lambda)$ with the section $\lb\to A(\lb)$ defines a continuous structure on $p:\mathcal{K}\to \Lambda$, and   $\mathbb{K}_0(\Lambda)$ becomes the corresponding space of continuous sections vanishing at infinity. In particular, the map $\lb\to\|A(\lb)\|$ is continuous and vanishing at infinity, for every $A\in \mathbb{K}_0(\Lambda)$.
\end{proposition}

\begin{proof}
To prove this, we have to specify a $\ast$-algebra $\AA$ of sections such that satisfies {\it i)} and {\it ii)} in definition \ref{contical}. We take $\AA$ to be the set of linear combinations of operators of the form $|\varphi\rangle \langle\psi|$\,, with $\varphi,\,\psi\in\Gamma_0^0(\Lambda)$.

Let $V$ be a Hilbert space and $\varphi_1,\cdots\varphi_{2n}$ be variable vectors. Assume that $\lim\langle\varphi_i,\varphi_j\rangle=\langle\psi_i,\psi_j\rangle$, where $\psi_1,\cdots,\psi_{2n}\in V$. Lemma 3.5.6 in \cite{D2} implies that $$\lim\left\||\varphi_1\rangle \langle\varphi_2|+\cdots+|\varphi_{2n-1}\rangle \langle\varphi_{2n}|\right\|=\left\||\psi_1\rangle \langle\psi_2|+\cdots+|\psi_{2n-1}\rangle \langle\psi_{2n}|\right\|.
$$ 
In particular, if $\langle\varphi_i,\varphi_j\rangle=\langle\psi_i,\psi_j\rangle,$ then $$\left\||\varphi_1\rangle \langle\varphi_2|+\cdots+|\varphi_{2n-1}\rangle \langle\varphi_{2n}|\right\|=\left\||\psi_1\rangle \langle\psi_2|+\cdots+|\psi_{2n-1}\rangle \langle\psi_{2n}|\right\|.$$ Therefore, the map $V^{2n}\ni(\psi_1,\cdots,\psi_{2n})\to \left\||\psi_1\rangle \langle\psi_2|+\cdots+|\psi_{2n-1}\rangle \langle\psi_{2n}|\right\|$ is actually a continuous function on the inner products $\langle\psi_i,\psi_j\rangle$. Thus, $\AA$ satisfies {\it{i)}} and clearly it also satisfies {\it{ii)}}. Moreover, since $\AA$ is closed under multiplication by elements of $C_0(\Lambda)$, the $\ast$-algebra $\AA$ is dense in the set of continuous sections vanishing at infinity of the corresponding $C^*$-bundle. By definition, the closure of $\AA$ is $\mathbb{K}_0(\Lambda)$, and this finishes the proof.

\end{proof}

In view of theorem \ref{desc}, proposition \ref{inter} and equation \eqref{Lei}, we define a smooth field of $C^*$-algebras as follows.

\begin{Definition}\label{DefC}
Let $\Lambda$ be a smooth manifold and $p: \A \rightarrow \Lambda$ be a field of $C^*$-algebras. A smooth structure on $\A\to\Lambda$ is given by specifying a $\ast$-algebra of sections $\AA^\infty$ closed under multiplication by elements of $C^\infty(\Lambda)$, and a map 
$\hat{\nabla}:\operatorname{Vect}(\Lambda)\times \AA^\infty\mapsto \AA^\infty$ such that for  $X,Y\in \operatorname{Vect}(\Lambda)$, $A,B\in\AA^\infty$ and $a\in C^{\infty}(\Lambda)$
\begin{enumerate}
\item[i)] $\hat\nabla_{X+Y}(A)=\hat\nabla_{X}(A)+\hat\nabla_{Y}(A)$ and  $\hat\nabla_{aX}(A)=a\nabla_{X}(A)$, 
\item[ii)] $\hat\nabla_{X}(aA)=X(a)A+a\hat\nabla_{X}(A)$ and $\hat\nabla_{X}(AB)=\hat\nabla_X(A)B+A\hat\nabla_{X}(B)$,
\item[iii)] $(\hat\nabla_{X}(A))^{*}=\hat\nabla_{\overline{X}}(A^{*})$,
\item[iv)] For each $m\in\mathbb{N}$ and $X_1,\cdots,X_m\in\operatorname{Vect}(\Lambda)$, the map $\lambda\mapsto \|\hat\nabla_{X_1}\cdots\hat\nabla_{X_m}A(\lb)\|$ is continuous,
\item[v)] $\AA^\infty(\lb)=\{A(\lb)\mid A\in\AA^\infty\}$ is dense in $\AA(\lb)$, for all $\lb\in\Lambda$.
\end{enumerate}
\end{Definition}

Let us return to the smooth fields of Hilbert spaces framework. A straightforward  
computation shows that, if $\varphi,\psi\in\Gamma^n(\Lambda)$ then $|\varphi\rangle \langle\psi|\in\AA^n_{lb}$ and the following natural identity holds
$$
\hat\nabla_X(|\varphi\rangle \langle\psi|)=|\nabla_X\varphi\rangle \langle\psi|+|\varphi\rangle \langle\nabla_X\psi|
$$
The latter fact and proposition \ref{exako} implies the following result. 

\begin{Corollary}\label{cor414} Let $H\to\Lambda$ be a smooth field of Hilbert spaces with connection $\nabla$ and let $\mathcal K\to\Lambda$ be the field of $C^*$algebras with fibers $\mathbb{K}(\Hi(\lb))$\,.
The space of sections $\mathbb{K}^\infty=\text{span}\{|\varphi\rangle \langle\psi|\mid \varphi,\psi\in \Gamma^{\infty}\} $ together with the connection $\hat\nabla=[\nabla,\cdot]$ makes $\mathcal K\to\Lambda$ a smooth field of $C^*$-algebras.
\end{Corollary}
  
Let $\varphi,\psi\in C_0(\Lambda,V)$, where $V$ is some Hilbert space. An straightforward computations show that
$$
\||\varphi\rangle\langle\psi|(\lb)-|\varphi\rangle\langle\psi|(\lb_0)\|\leq \|\varphi(\lb)-\varphi(\lb_0)\|\|\psi(\lb)\|+\|\psi(\lb)-\psi(\lb_0)\|\|\varphi(\lb_0)\|.
$$
Therefore, any finite rank operator on $\Gamma_0$ belongs to $C_0(\Lambda,\mathbb{K}(V))$, where $\mathbb{K}(V)$ is endowed with the norm operator topology. The following result holds within the more general framework of continuous fields of Hilbert spaces. We use the notion of continuous trivialization of Dixmier and Douady (definitions 2 and 3 in \cite{D3}).

\begin{proposition}
Let $H\to\Lambda$ be a continuous field of Hilbert spaces and $\Gamma_0$ be the corresponding space of continuous sections vanishing at infinity. Also, let $T:H\to V$ be a continuous trivialization. The map $\hat T:\mathbb K(\Gamma_0)\to C_0(\Lambda,\mathbb{K}(V))$ given by $\hat T(A)(\lb)=T_\lb A(\lb)T^*_\lb$ is an isomorphism of $C^*$-algebras.
\end{proposition}
\begin{proof}
Notice that $\hat T(|\varphi\rangle\langle\psi|)=|T(\varphi)\rangle\langle T(\psi)|$. Therefore $\hat T(A)\in C_0(\Lambda,\mathbb{K}(V))$, for any finite rank operator $A$ on $\Gamma_0$. Since $\hat T$ is an isometry, $\hat T$ is well-defined. Moreover, since $T(\Gamma_0)=C_0(\Lambda,V)$ (fullness is included in the definition of continuous trivialization), $\hat T(\mathbb K(\Gamma_0))(\lb)=\mathbb K(V)$. Therefore, proposition C.24 in \cite{Wi} implies that $\hat T(\mathbb K(\Gamma_0(\Lambda)))$ is dense in $C_0(\Lambda,\mathbb{K}(V))$ and so $\hat T$ is surjective.
\end{proof}

\appendix
\section{Measurable and Continuous fields of Hilbert Spaces}\label{MyC}
In this appendix we recall the definitions of measurable and continuous fields of Hilbert spaces over a space $\Lambda$. Each of those notions admit at least three different but equivalent ways to introduce them. Let us discuss some of their main characteristics in order to understand their relation.

\medskip
The first (historic) definition of a measurable field of Hilbert spaces describes the measurable structure in terms of a space of sections. In fact, the latter approach is used throughout this article. 
\begin{Definition}\label{mfhs}
Let $\Lambda$ be a measurable space. A measurable structure on a field of Hilbert spaces $p:H\to\Lambda$ is given by specifying a linear space of sections $\Gamma$ possessing the following properties:
\begin{description}
  \item (i) The function $\lb\to ||\varphi(\lb)||_\lb$ is measurable, for every $\varphi\in \Gamma$. 
  \item (ii) If $\psi$ is a section such that, for every $\varphi\in\Gamma$, the function
              $\lb\to \langle \varphi(\lb),\psi(\lb)\rangle_\lb$ is measurable, then $\psi\in\Gamma$.
  \item (iii) There exists a sequence $(\varphi_j)$ of elements of  $\Gamma$ such that $\mbox{span}\{\varphi_j(\lb)\mid j \in \mathbb{N}\}$ is dense in $\Hi(\lb)$, for every $\lb\in\Lambda$.
\end{description}
The elements of $\Gamma$ are called measurable sections. The sequence $(\varphi_1,\varphi_2,\cdots)$ is called a base of measurability or a measurable field of bases. The set
$\{\lb\in\Lambda \mid \Hi(\lb)\neq 0\}$ is called the support of the field.
\end{Definition}
The second possible way to introduce measurability into fields of Hilbert spaces is given by the notion of a measurable Hilbert bundle (definition \ref{mhb} below or definition 2.4.8.\ in \cite{We}). Let us explain how we meet that notion beginning with a measurable field of Hilbert spaces.  

Starting from the sequence $(\varphi_j)$ we can obtain a local field of measurable orthonormal bases. More precisely, there is a measurable partition $\Lambda=\bigcup_{n=1}^{n=\infty} \Lambda_n$ and a sequence $\phi_j:\Lambda\to H$ such that, $\text{dim }\Hi(\lambda)=n$ for all $\lambda\in\Lambda_n$,  $\phi_j(\lambda)=0$ for all $j>\text{dim }\Hi(\lambda)$ and $\{\phi_j(\lambda)\mid 1\leq j\leq \text{dim }\Hi(\lambda)\}$ is a orthogonal basis of $\Hi(\lambda)$. Up to certain technicalities, the main idea is just to apply Gram-Schmidt orthogonalization process. Details can be found in part II chapter 1, lemma 2.1 and proposition 4.1 in \cite{D}. Since choosing a pointwise orthonormal basis on the fiber $\Hi(\lambda)$ is equivalent to fix an unitary operator from $\Hi(\lambda)$ to $l^{2}(\mathbb{N})$, then the partition induces maps $T_n:H|_{\Lambda_n}\to \Hi_n$ such that $\Hi_n$ is a Hilbert space and $T_n|_{\Hi(\lambda)}$ is unitary for all $\lambda\in\Lambda_n$. Moreover, the map  $T:H\to \bigsqcup \Hi_n$ defined by $T(x)=T_n(x)$, for every $x\in H$ with $p(x)\in\Lambda_n$, can be interpreted as a sort of measurable trivialization. Indeed, $\varphi$ is measurable if and only if $T(\varphi)$ is weakly measurable. Therefore $H$ is equivalent to the bundle $\bigcup \Lambda_n \times \Hi_n$, and we have met the following definition of measurable Hilbert bundle.
\begin{Definition}\label{mhb}
Let $\Lambda$ be a measurable space. A (separable) measurable Hilbert bundle (MHB) over $\Lambda$ is a disjoint union $H=\bigcup \Lambda_n \times \Hi_n$ where $\{\Lambda_n\}$ is a measurable partition of $\Lambda$ and $\Hi_n$ is a Hilbert space of dimension $n$, with $0< n\leq \infty$.
\end{Definition}
Conversely, if $T$ is a trivialization $T:H\to \bigsqcup \Hi_n$ as before, we can define the space $\Gamma$ as the set of all the sections $\varphi$ of $H$ such that $T\varphi$ is weakly measurable. Clearly, $\Gamma$ satisfies \textit{(i)}, \textit{(ii)} and \textit{(iii)} of definition \ref{mfhs}.\\
A much richer structure is obtained if a fixed measure $\eta$ on $\Lambda$ is considered. Notice that the notion of measurable Hilbert bundle in \cite{We} requires $\Lambda$ to be a $\sigma$-finite measure space. Moreover, the third way to describe measurable fields of Hilbert spaces also requires the latter additional assumption and it is given within the context of Hilbert modules over $L^{\infty}(\Lambda,\eta)$. If $\mathcal{E}$ is such Hilbert module, we can define the dual module $\mathcal{E}^{*}$ in a canonical way, i.e. taking the bounded linear maps from the module to $L^{\infty}(\Lambda,\eta)$. When every element $f\in \mathcal{E}^{*}$ admits a representation of the form $f(\phi)=\langle \phi,\psi\rangle$ for some $\psi\in \mathcal{E}$, the module is called self-dual. If $\Lambda$ is a $\sigma$-finite measure space and $H\to \Lambda$ is a measurable Hilbert bundle,  it is not difficult to show that  the space of weakly measurable essentially bounded sections of a measurable bundle is a self-dual Hilbert $L^{\infty}(\Lambda)$-module (propositions 9.2.3 in \cite{We}). Conversely, every self-dual weakly separable module is obtained in the latter way (theorem 9.2.4 in \cite{We}). In summary, there are two equivalent ways to introduce measurability for fields of Hilbert spaces over a measurable space $\Lambda$.
\begin{enumerate}
    \item Measurable fields of Hilbert spaces.
    \item Measurable Hilbert bundles.
\end{enumerate}
If in addition we assume $(\Lambda,\eta)$ is a $\sigma$-finite measure space the latter notions are equivalent to the following.
\begin{enumerate}
    \item[3.] Self-dual weakly separable Hilbert $L^{\infty}(\Lambda,\eta)$-modules.
\end{enumerate}
One of the main reasons to consider measurable fields of Hilbert spaces is that it allow us to introduce direct integrals (we require them in subsection \ref{IMS}).
\begin{Definition}
Let $(\Lambda,\eta)$ be a measure space and $H\to \Lambda$ be a measurable field of Hilbert spaces. Up to quotient by the space of measurable sections vanishing almost everywhere, the direct integral $\int^\bigoplus_\Lambda \Hi(\lb)\de\eta(\lb)$ is the space of square
integrable measurable vector fields, i.e.
$$
\int^\bigoplus_\Lambda \Hi(\lb)\de\eta(\lb):=\left\{\varphi\in\mathcal S \left| \right.   \int_\Lambda|| \varphi(\lb)||^2_\lb\de\eta(\lb)<\infty \right \}
$$
\end{Definition}
It is well-known (for instance, see part II, chapter I, proposition 5 in \cite{D}) that
$\int^\bigoplus_\Lambda \Hi(\lb)\de\eta(\lb)$ is a Hilbert space with the inner product
$$
\langle\varphi,\psi\rangle:=\int_\Lambda \langle\varphi(\lb),\psi(\lb)\rangle_\lb\de\eta(\lb).
$$
In the topological framework there are also three equivalent ways to introduce continuity for fields of Hilbert spaces.\ Assume that $\Lambda$ is a Hausdorff locally compact space. The first historical definition given by Godement \cite{God} of a continuous field of Hilbert spaces is the following.
\begin{Definition}\label{CField}
Let $\Lambda$ be a locally compact space and $H \to \Lambda$ a field of Hilbert spaces. A continuous structure is given by specifying a linear space of sections $\Gamma$ possesing the following properties
\begin{description}
    \item (i) The set $\{\varphi(\lambda) \mid \varphi \in \Gamma \}$ is dense in $\Hi(\lambda)$, for every $\lambda\in\Lambda$ .
    \item (ii) The function $\lambda\to\|\varphi(\lambda)\|$ is continuous, for every $\varphi\in\Gamma$.
\end{description}
\end{Definition}

Latter Dixmier and Douady in \cite{D2,D3} added the following condition.
\begin{description}
    \item \textit{(iii)} If $\varphi$ is a section and for every $\lambda_0\in \Lambda$ and every $\varepsilon >0$ there exists $\varphi'\in\Gamma$ such that $\|\varphi(\lambda)-\varphi'(\lambda)\|\leq \varepsilon$ for every $\lambda$ in some neighborhood (depending on $\varepsilon$) of $\lambda_0$, then $\varphi\in\Gamma$.
\end{description}
They proved that if $\Gamma$ satisfies conditions \textit{(i)} and \textit{(ii)} then there exists a unique space of sections $\tilde \Gamma$ satisfying \textit{1)}, \textit{2)}  and \textit{3)} such that $\Gamma \subseteq \tilde\Gamma$ (proposition 3 in \cite{D3} or proposition 10.2.3.\ in \cite{D2}). Moreover, they also show that if $\Gamma$ satisfies \textit{1)}, \textit{2)} and \textit{3)} then the following properties hold. 
\begin{enumerate}
    \item[\textit{iv)}] $\Gamma$ is a $C(\Lambda)$-module.
    \item[\textit{v)}] The set $\{\varphi(\lambda) \mid \varphi \in \Gamma \}$ equals $\Hi(\lambda)$, for every $\lambda\in\Lambda$.
\end{enumerate}
Note that if $\Gamma$ satisfies  \textit{1)}, \textit{2)} and \textit{4)} then $\Gamma$ is dense in $\tilde\Gamma$ with respect to the local uniform convergence topology. 

Now we consider the definition of a continuous Hilbert bundle. The first definition of this structure can be found in \cite{TkJ}. We follow \cite{We} as reference as we did in the measurable case.
\begin{Definition}
Let $\Lambda$ be a compact Hausdorff space. A covering space of $\Lambda$ is a topological space $H$ together with a continuous open surjection $p : H \to \Lambda$. A continuous Hilbert bundle over $\Lambda$ is then a covering space $H$ such that $\mathcal{H}(\lambda) = p^{-1}(\lambda)$ is equipped with a Hilbert space structure for each $\lambda\in\Lambda$, and satisfying the following conditions:
\begin{description}
    \item (i) the map $x \to \|x\|$ is continuous from $H$ to $\R$;
    \item (ii) the map $(x_1, x_2) \to x_1 + x_2$ is continuous from $H\times H$ to $H$; 
    \item (iii) the map $x \to ax$ is continuous from $H$ to $H$ for every $a \in \C$; and
    \item (iv) for any neighborhood $\mathcal{O}$ of the origin of $\Hi(\lambda)$ in $H$ there exists a neighborhood $\mathcal{O}'$ of $\lambda$ in $\Lambda$ and an $\varepsilon > 0$ such that 
    $$\{ x\in H : p(x)\in \mathcal{O}' \mbox{ and } \|x\| < \varepsilon\} \subset \mathcal{O}.$$
\end{description}
 A section of $H$ is a function $\varphi: \Lambda \to H$ such that $\varphi(\lambda) \in \mathcal{H}(\lambda)$, for all $\lambda \in \Lambda$. The set of all continuous sections of $H$ is denoted $\Gamma(H)$.
\end{Definition}

Clearly, if $H\to\Lambda$ is a continuous Hilbert bundle, then $\Gamma(H)$ satisfies the conditions {\it (i), (ii)} and {\it (iii) } of definition \ref{CField}. Conversely, if $\Gamma$ satisfies {\it (i)} and {\it (ii)} of definition \ref{CField} then there is a (final) topology on $H$ such that $H\to\Lambda$ is a continuous Hilbert bundle and $\Gamma\subset\Gamma(H)$ (in fact $\tilde\Gamma=\Gamma(H)$). 

\medskip
The third way to describe continuous field of Hilbert spaces is through the notion of Hilbert $C_0(\Lambda)$-module. In fact, the space of continuous sections vanishing at infinity $\Gamma_0(H)$ of a continuous Hilbert bundle is clearly a Hilbert $C_0(\Lambda)$-module. Conversely, every Hilbert $C_0(\Lambda)$-module $\Gamma_0$ corresponds to the space of continuous sections vanishing at infinity of a Hilbert bundle. Indeed, for $\lambda\in\Lambda$, let  $I_{\lambda}=\{f\in C_0(\Lambda)\mid f(\lambda)=0\}$\,. It is easy to show that the quotient $\mathcal{H}(\lambda)=\faktor{\Gamma_0}{I_{\lambda}\Gamma_0}$ is a Hilbert space. Therefore, $\Gamma_0$ can be regarded as a space of sections by defining $\varphi(\lambda)=\pi_\lambda(\varphi)$ for every $\varphi\in\Gamma_0$, where $\pi_\lambda:\Gamma_0\to \faktor{\Gamma_0}{I_{\lambda}\Gamma_0}$ is the canonical projection and $\varphi\in\Gamma_0$. Moreover, the map $\lambda\to\|\varphi(\lambda)\|$ is clearly continuous, then $\Gamma_0$ defines a continuous field of Hilbert spaces (see definition \ref{CField}). The latter equivalence was first noticed in \cite{TkJ}. If we want to pass directly from Hilbert $C_0(\Lambda)$-modules to continuous Hilbert bundles over $\Lambda$, we can apply  proposition 9.15 in \cite{We} in the compact case. In the general case and use a one-point compactification, as explained in \cite{Ph}. 

Summarizing, there are three equivalent ways to introduce the notion of continuity for fields of Hilbert spaces.
\begin{enumerate}
    \item Continuous fields of Hilbert spaces.
    \item Continuous Hilbert bundles.
    \item Hilbert $C_0(\Lambda)$-modules.
    
\end{enumerate}

The reader should notice the similarities between these three descriptions and the analogous three descriptions given for the measurable framework. 

In subsection \ref{LB} a similar discussion can be found for the case of continuous fields of $C^{*}$-algebras.


\end{document}